\documentclass[reqno,11pt]{amsart} 
\usepackage{amsmath}
\usepackage{amssymb}
\usepackage{amsthm}
\usepackage{hyperref}
\usepackage{amsfonts,wasysym,mathrsfs}
\usepackage{xcolor}
\usepackage{graphics}

\newcommand\NoBlackBoxes{\global\overfullrule0pt}
\NoBlackBoxes
\theoremstyle{plain}

\newtheorem{theorem}{Theorem}[section]
\newtheorem{lemma}[theorem]{Lemma}
\newtheorem{corollary}[theorem]{Corollary}

\newtheorem{definition}[theorem]{Definition}

\numberwithin{equation}{section}

\begin{document}

\title{On Gilles Pisier's approach to Gaussian concentration, \\
isoperimetry, and Poincar\'e-type inequalities}

\author{Sergey G. Bobkov$^{1}$}
\thanks{1) 
School of Mathematics, University of Minnesota, Minneapolis, MN, USA
}

\author{Bruno Volzone$^{2}$}
\thanks{2) Dipartimento di Scienze Economiche, Giuridiche, Informatiche e Motorie, Universit\`a degli Studi di Napoli ``Parthenope'',  Centro Direzionale Isola C4, 80143 Napoli, Italy.
}


\begin{abstract}
We discuss a natural extension of Gilles Pisier's approach to the study of measure
concentration, isoperimetry, and Poincar\'e-type inequalities. This approach allows one
to explore counterparts of various results about Gaussian measures in the class of
rotationally invariant probability distributions on Euclidean spaces,
including multidimensional Cauchy measures.
\end{abstract} 
 
\maketitle
\markboth{Sergey G. Bobkov and Bruno Volzone}{Gilles Pisier approach}

\def\theequation{\thesection.\arabic{equation}}
\def\E{{\mathbb E}}
\def\R{{\mathbb R}}
\def\C{{\mathbb C}}
\def\P{{\mathbb P}}
\def\Z{{\mathbb Z}}
\def\S{{\mathbb S}}
\def\I{{\mathbb I}}
\def\T{{\mathbb T}}

\def\s{{\mathbb s}}

\def\G{\Gamma}

\def\Ent{{\rm Ent}}
\def\var{{\rm Var}}
\def\Var{{\rm Var}}
\def\cov{{\rm cov}}
\def\V{{\rm V}}

\def\H{{\rm H}}
\def\Im{{\rm Im}}
\def\Tr{{\rm Tr}}
\def\s{{\mathfrak s}}
\def\A{{\mathfrak A}}
\def\m{{\mathfrak m}}

\def\k{{\kappa}}
\def\M{{\cal M}}
\def\Var{{\rm Var}}
\def\Ent{{\rm Ent}}
\def\O{{\rm Osc}_\mu}

\def\ep{\varepsilon}
\def\phi{\varphi}
\def\vp{\varphi}
\def\F{{\cal F}}

\def\be{\begin{equation}}
\def\en{\end{equation}}
\def\bee{\begin{eqnarray*}}
\def\ene{\end{eqnarray*}}

\thispagestyle{empty}

\tableofcontents














\section{{\bf G. Pisier's approach and its consequences}}\label{PISIER}
\setcounter{equation}{0}

\vskip2mm
\noindent
Let $\gamma_n$ denote the standard Gaussian measure on $\R^n$,
thus with density 
$$
\frac{d\gamma_n(x)}{dx} = (2\pi)^{-\frac{n}{2}}\, e^{-|x|^2/2}, \quad x \in \R^n,
$$
with respect to the Lebesgue measure, where $|\cdot|$ stands for the canonical Euclidean norm. 
In the mid 1980's G. Pisier \cite{P} proposed the following remarkable family of integro-differential
inequalities involving this measure.

\begin{theorem}\label{gaussianineqth} Let $\Psi:\R \rightarrow \R$ be a convex function. For any smooth
function $f$ on $\R^n$ with gradient $\nabla f$,
\begin{equation}
\int_{\R^n} \int_{\R^n} \Psi(f(y) - f(x))\,d\gamma_n(x)\,d\gamma_n(y)
\leq \int_{\R^n} \int_{\R^n} 
\Psi\Big(\frac{\pi}{2} \left<\nabla f(x),y\right>\Big)\,d\gamma_n(x)\,d\gamma_n(y).\label{Pisier}
\end{equation}
In particular, if $f$ has $\gamma_n$-mean zero, then
\begin{equation}
\int_{\R^n} \Psi(f)\,d\gamma_n \leq \int_{\R^n} \int_{\R^n} 
\Psi\Big(\frac{\pi}{2} \left<\nabla f(x),y\right>\Big)\,d\gamma_n(x)\,d\gamma_n(y).\label{Pisier2}
\end{equation}
\end{theorem}
\vskip2mm
This result was actually stated in a more general setting of Gaussian measures on 
locally convex spaces, which is readily suggested by the dimension-free character of \eqref{Pisier}.

What is also surprising, this inequality admits a rather simple proof.
For reader's convenience, we include it below with simplications
due to B. Maurey (according to \cite{P}). 

\vskip5mm
{\bf Proof.} Given independent random vectors $X$ and $Y$ in $\R^n$ with distribution
$\gamma_n$, put $X(t) = X \cos t + Y \sin t$ for $0 \leq t \leq \frac{\pi}{2}$. 
By the Leibniz formula,
$$
\Delta \equiv f(Y) - f(X) = \int_0^{\pi/2} \frac{d}{dt} f(X(t))\,dt =\int_0^{\pi/2}
\left<\nabla f(X(t)),X'(t)\right>dt,
$$
where $X'(t) = -X \sin t + Y \cos t$. Hence, by Jensen's inequality,
$$
\Psi(\Delta) \leq \frac{2}{\pi} \int_0^{\pi/2} \Psi\Big(\frac{\pi}{2}
\left<\nabla f(X(t)),X'(t)\right>\Big)\,dt.
$$
Taking the expectation, we get
\be
\E\,\Psi(\Delta) \leq \frac{2}{\pi} \int_0^{\pi/2} \E\, \Psi\Big(\frac{\pi}{2}\label{Pisier3}
\left<\nabla f(X(t)),X'(t)\right>\Big)\,dt.
\en

Now, a crucial point is that, since  the Gaussian measure 
$\gamma_{2n} = \gamma_n \otimes \gamma_n$ is rotationally invariant on 
$\R^{2n}$, the couple $(X(t),X'(t))$ represents an independent copy of $(X,Y)$.
In particular, the second expectation in \eqref{Pisier3} in the integrand does not depend on $t$.
\qed

\vskip2mm
The relations \eqref{Pisier}-\eqref{Pisier2} have several important consequences.
Choosing $\Psi(r) = |r|^p$ with $p \geq 1$, we get a Poincar\'e-type inequality 
for $L^p$-norms in Gauss space $(\R^n,\gamma_n)$. Namely,
\be
\int_{\R^n} \int_{\R^n} |f(x) - f(y)|^p\,d\gamma_n(x)\,d\gamma_n(x) \, \leq \,\label{Gaussdoubineq}
c_p \int_{\R^n} |\nabla f|^p \, d\gamma_n,
\en
and, by Jensen's inequality,
\be
\int_{\R^n} \Big|f - \int_{\R^n} f\, d\gamma_n\Big|^p\,d\gamma_n \, \leq \, \label{Gaussclassineq}
c_p \int_{\R^n} |\nabla f|^p \, d\gamma_n
\en
with
$$
c_p = \Big(\frac{\pi}{2}\Big)^p\,\E\,|\xi|^p = 
\Big(\frac{\pi}{2}\Big)^p\ \frac{2^{\frac{p}{2}}}{\sqrt{\pi}}\,
\Gamma\Big(\frac{p+1}{2}\Big),
$$
where $\xi$ is a normal random variable with distribution $\gamma_1$.
Such relations are known to hold for general measures as a consequence of 
the particular case $p=1$ when the latter is true -- however, with constants 
growing like $(cp)^p$ for large $p$ (where $c>0$ is an absolute constant). 
In the present formulation, 
$c_p$ is of the order $(cp)^{p/2}$ which is asymptotically correct. 

While the best constants $c_p$ are unknown except for $p = 1,2$, the constant 
$c_1 = \sqrt{\frac{\pi}{2}}$ is sharp. In this case, \eqref{Gaussdoubineq}-\eqref{Gaussclassineq} are 
equivalent to the Cheeger-type isoperimetric inequality
$$
\gamma_n^+(A) \geq \sqrt{\frac{2}{\pi}}
\min\{\gamma_n(A),1-\gamma_n(A)\} \quad (A \subset \R^n \, {\rm Borel}),
$$
where $\gamma_n^+$ denotes the Gaussian perimeter. This relation becomes 
an equality for any half-space $A$ whose boundary contains the origin.

Another choice $\Psi(r) = e^r$ leads in \eqref{Pisier2} to the exponential bound
\be
\E\,e^{f(X)} \leq \E\,e^{\frac{\pi^2}{8}\, |\nabla f(X)|^2} \quad\label{exponentialboundgaus}
(\E f(X) = 0),
\en
which is one of the ways to express the Gaussian dimension-free concentration penomenon
(later in \cite{B-G}, the constant $\frac{\pi^2}{8}$ was removed from the exponent
by applying the logarithmic Sobolev inequality).
In particular, if the function $f$ is 1-Lipschitz, \eqref{exponentialboundgaus} implies that the random variable 
$f(X)$ is subgaussian, i.e. $\E \exp\{c f(X)^2\} \leq 2$ for some
constant $c>0$.

\vskip5mm
\section{{\bf Extension to general measures. Cauchy distributions}}\label{EXTENSIONS GENERAL MEASURES}
\setcounter{equation}{0}

\vskip2mm
\noindent
The integration in \eqref{Pisier} is carried out with respect to the Gaussian measure
$\gamma_{2n} = \gamma_n \otimes \gamma_n$ on $\R^n \times \R^n$.
This inspires a natural idea to explore the applicability of Pisier's approach
to general positive measures $\mu$ on $\R^n \times \R^n$. Recently, 
motivated by the Enflo problem, Ivanisvili, van Handel and Volberg \cite{I-H-V} 
have found a dimension-free analogue of Pisier's inequality \eqref{Pisier2} on the discrete 
cube $\{-1,1\}^n$ for the Bernoulli distribution and for functions $f$ with 
values in an arbitrary Banach space. Alternatively, we directly follow Pisier's
argument and introduce the following transform, without
requiring that the measures have a product structure. 

\vskip5mm
\begin{definition}\label{sphericalcup}Denote by $\mu_t$ the image
of $\mu$ under the orthogonal linear transformation
$$
U_t(x,y) = (x \cos t + y \sin t,-x \sin t + y \cos t), \quad x,y \in \R^n,
$$
from $\R^n \times \R^n$ into itself
(i.e. $U_t$ pushes forward $\mu$ onto $\mu_t$). We call the measure on
$\R^n \times \R^n$
\be
\hat \mu = \frac{2}{\pi} \int_0^{\frac{\pi}{2}} \mu_t\,dt\label{sphericalcup}
\en
a spherical cup for $\mu$.
\end{definition}

\vskip2mm
For example, a spherical cup for $\gamma_{2n}$ is $\gamma_{2n}$ itself.
More generally, $\hat \mu = \mu$ as long as $\mu$ is rotationally invariant, 
although such measures lose the product structure in the non-Gaussian case. 
With this in mind, Theorem 1.1 has the following extension.
Without loss of generality (and in order to avoid integrability questions), 
we assume that the convex function $\Psi$ is non-negative.

\vskip5mm
\begin{theorem}\label{maintheorem}
Given a positive measure $\mu$ on $\R^n \times \R^n$
with a spherical cup $\hat \mu$, for any smooth function $f$ on $\R^n$, 
\begin{equation}\label{maininequal}
\int_{\R^n \times \R^n} \Psi(f(y) - f(x))\,d\mu(x,y) \leq
\int_{\R^n \times \R^n} \Psi\Big(\frac{\pi}{2} \left<\nabla f(u),v\right>\Big)\,
d\hat \mu(u,v).
\end{equation}
\end{theorem}

\vskip2mm
In particular, this holds true with $\hat \mu = \mu$ if $\mu$ is rotationally invariant.
We discuss the transform $\mu \rightarrow \hat \mu$ and a similar
transformation for densities in Section \ref{SPHERICAL CUPS}.
Under certain isotropy-type conditions, the right integral in \eqref{maininequal} can be simplified
for the quadratic function $\Psi(r) = r^2$, which leads to weighted Poincar\'e-type
inequalities (Section \ref{ISOTROPIC MEASURES}). Moreover, $L^p$-Poincar\'e-type inequalities with arbitrary $p \geq 1$
can be obtained from \eqref{maininequal} for the class of rotationally invariant measures.
The example of the uniform distribution on the sphere is discussed in Section \ref{SPHERICALLY INVARIANT MEASURES}.

We also illustrate this kind of applications on the example of Cauchy (also called Student)
probability distributions. Recall that the $n$-dimensional Cauchy measure 
$\m_{n,\alpha}$ on $\R^n$ of order $\alpha > \frac{n}{2}$ 
has density proportional to $(1 + |x|^2)^{-\alpha}$, $x \in \R^n$.
We will remind of basic properties and identities related to this class
of probability distributions in Section \ref{BACKGROUND CAUCHY}. Here let us only mention that, similarly to
the Gaussian case, if the couple of random vectors $(X,Y)$ in $\R^{2n}$ has 
a Cauchy distribution of order $\alpha > n$, then both $X$ and $Y$ have 
a Cauchy distribution on $\R^n$ of order $\alpha - \frac{n}{2}$. 
As a Cauchy-type analog of \eqref{Gaussdoubineq}, we have:

\begin{theorem}\label{Poinccauchytheo} Let $\alpha>n+\frac{1}{2}$ and $1 \leq p < 2\,(\alpha - n)$.
For any smooth function $f$ on $\R^n$,
\begin{equation}\label{Poinccauchy}
\int_{\R^n \times \R^n} |f(x) - f(y)|^p\,d\m_{2n,\alpha}(x,y) \, \leq \,
C \Big(\frac{\pi}{2}\Big)^p
\int_{\R^n} |\nabla f(x)|^p \, d\m_{n,\beta}(x),
\end{equation}
where $\beta = \alpha - \frac{n+p}{2}$, and where the constant depends
on $(n,p,\alpha)$ and is given by
$$
C = \frac{1}{\sqrt{\pi}}\,
\frac{\Gamma(\frac{p+1}{2})\Gamma(\alpha - n - \frac{p}{2})}{\Gamma(\alpha - n)}.
$$
\end{theorem}

\vskip2mm
In particular, for $p=2$ and with $\beta = \alpha - \frac{n}{2} - 1$, \eqref{Poinccauchy} takes the form
\begin{equation}
\int_{\R^n \times \R^n} |f(x) - f(y)|^2\,d\m_{2n,\alpha}(x,y) \, \leq \,
\frac{\pi^2}{8\, (\alpha - n - 1)} \int_{\R^n} |\nabla f|^2 \, d\m_{n,\beta}.\label{Poinccauchyp=2}
\end{equation}

The right integral in \eqref{Poinccauchy} may be written over the measure $\m_{n,\alpha}$
with weight $(1 + |x|^2)^{\frac{n+p}{2}}$ using a different constant in place of $C$.
Hence, this relation may be viewed as a weighted Poincar\'e-type inequality for
$L^p$-norms.

If $\alpha$ is sufficiently large, the constants in \eqref{Poinccauchy}-\eqref{Poinccauchyp=2} do not exceed 
a multiple of $1/\alpha$, which is a correct growth rate. Moreover, in terms of 
the image $\widetilde \m_{2n,\alpha}$ of $\m_{2n,\alpha}$ under the linear map 
$(x,y) \rightarrow \sqrt{2\alpha}\,(x,y)$ and the image $\widetilde \m_{n,\beta}$ 
of $\m_{n,\beta}$ under the map 
$x \rightarrow \sqrt{2\alpha}\,x$, that is the measures with the densities
\[
d\widetilde \m_{2n,\alpha}(x,y)=\frac{(2\alpha)^{-n}}{c_{2n,\alpha}}\left(1+\frac{|x|^{2}+|y|^{2}}{2\alpha}\right)^{-\alpha}dx\, dy,
\]
\[
d\widetilde \m_{n,\beta}(x)=\frac{(2\alpha)^{-n/2}}{c_{n,\beta}}\left(1+\frac{|x|^{2}}{2\alpha}\right)^{-\beta}dx
\]
for suitable normalization constants $c_{2n,\alpha}$, $c_{n,\beta}$ (see \eqref{normalizconst}), inequality \eqref{Poinccauchy} becomes
\begin{equation}\label{rescaledCauchy}
\int_{\R^n \times \R^n} |f(x) - f(y)|^p\,d\widetilde \m_{2n,\alpha}(x,y) \, \leq \,
C\,(2\alpha)^{\frac{p}{2}}\, \Big(\frac{\pi}{2}\Big)^p
\int_{\R^n} |\nabla f|^p \, d\widetilde \m_{n,\beta}.
\end{equation}
In the limit as $\alpha \rightarrow \infty$, 
$\widetilde \m_{2n,\alpha} \rightarrow \gamma_{2n}$ and
$\widetilde \m_{n,\beta} \rightarrow \gamma_n$ in the weak topology
(and actually in total variation distance),
while $C\,(2\alpha)^{\frac{p}{2}} \rightarrow 2^{p/2}\pi^{-1/2}\,\Gamma((p+1)/2)$. Hence, \eqref{rescaledCauchy}
recovers Pisier's Poincare-type inequality \eqref{Gaussdoubineq}. In this sense Cauchy's distributions 
after a proper normalization represent a pre-Gaussian model, which may be used
to recover various relations for the Gaussian measure.

On the other hand, this class of probability distributions is of independent interest
and has been intensively studied in the literature. In particular, inequalities such as
\eqref{Poinccauchyp=2}, i.e. \eqref{Gaussdoubineq} for $p=2$, were considered with respect to the product measure
$\m_{n,\alpha} \otimes \m_{n,\alpha}$ on the left-hand side, and with weight 
$1 + |x|^2$ over $\m_{n,\alpha}$ on the right-hand side, cf. e.g. \cite{B-L1}, \cite{B-L2}. 
For comparison, several results in this direction will be discussed at the end 
of this paper. One should mention here that the $2n$-dimensional Cauchy
measure $\m_{2n,\alpha}$ is rather close to the product 
$\m_{n,\alpha} \otimes \m_{n,\alpha}$ for the growing parameter $\alpha$.
In particular,
\begin{equation}\label{boundofcauchy}
\m_{2n,\alpha} \geq d\,\m_{n,\alpha} \otimes \m_{n,\alpha}, \quad
d = d_{n,\alpha} = 
\frac{\Gamma(\alpha - \frac{n}{2})^2}{\Gamma(\alpha - n)\, \Gamma(\alpha)}.
\end{equation} 
For example, $d \geq \frac{1}{2}$ for $\alpha \geq n^2$.\\

\noindent Of a special interest in \eqref{Poinccauchy} is the case $p=1$, when this inequality yields
\begin{equation}
\int_{\R^n \times \R^n} |f(x) - f(y)|\,d\m_{2n,\alpha}(x,y) \, \leq \,
\frac{\sqrt{\pi}}{\sqrt{\alpha-n}} \int_{\R^n} |\nabla f| \, d\m_{n,\beta}\label{Cauchyp=1}
\end{equation}
with $\alpha \geq n+1$ and $\beta = \alpha - \frac{n+1}{2}$.
Being applied to (nearly) indicator functions, \eqref{Cauchyp=1} together with \eqref{boundofcauchy} lead 
to the Cheeger-type isoperimetric-type inequality for the Cauchy measures. 
Let us recall that
the $\nu$-perimeter, or the outer Minkowski content for a (Borel) probability
measure $\nu$ on $\R^n$ is defined for all Borel sets $A$ in $\R^n$ by
\begin{equation}\label{perimeter}
\nu^+(A) = \liminf_{\ep \rightarrow 0}\, \frac{\nu(A_\ep) - \nu(A)}{\ep},
\end{equation}
where $A_\ep$  $(\ep>0)$ denotes an open $\ep$-neighborhood of $A$ for the
Euclidean distance.

\vskip5mm
\begin{corollary}\label{corisopercauchymeasure} Let $\beta \geq \frac{n+1}{2}$ and 
$\beta^* = \beta + \frac{n+1}{2}$. For any Borel set $A$ in $\R^n$,
\begin{equation}
\m_{n,\beta}^+(A) \geq c\, \m_{n,\beta^*}(A) (1 - \m_{n,\beta^*}(A)),\label{isopineqcauc}
\end{equation}
where one may take $c = \frac{d}{\sqrt{\pi}}\sqrt{2\beta}$, so that
$c \geq \frac{1}{\sqrt{2\pi}}\sqrt{\beta}$, if $\beta \geq n^2$.
\end{corollary}

\vskip5mm
Usually, the isoperimetric problem is aimed to estimate the perimeter
$\nu^+(A)$ via the ``size" $\nu(A)$. Here, however, for the lower
bound on $\m_{n,\beta}^+(A)$ we use a different Cauchy measure
$\m_{n,\beta^*}$ with a good isoperimetric profile function.
In particular, if $\m_{n,\beta^*}(A) = \frac{1}{2}$, the $\m_{n,\beta}$-perimeter
of $A$ is bounded from below by a large quantity for large values of $\beta$. This is
consistent with the isoperimetric inequality \eqref{exponentialboundgaus} for the Gaussian measure.

Finally, let us emphasize a concentration aspect behind the Poincar\'e-type 
inequality \eqref{Poinccauchy} generalizing the dimension free 
concentration phenomenon for the Gaussian measure. The latter can be stated
as the property that, for any function $f$ on $\R^n$ with Lipschitz semi-norm 
$\|f\|_{\rm Lip} \leq 1$, the function $f(x) - f(y)$ is subgaussian under
$\gamma_{2n}$, i.e.
$$
\gamma_{2n}\big\{(x,y)\in \R^n \times \R^n: 
|f(x) - f(y)| \geq t\big\} \leq 2\,e^{-ct^2}, \quad t>0,
$$
with some absolute constant $c>0$. A similar assertion continues to
hold for the Cauchy measures after the natural normalization of the
space variable.

\vskip5mm
\begin{corollary}\label{corlargedev}
If $\alpha\geq n+3/2$, for any function $f$ on 
$\R^n$ with $\|f\|_{\rm Lip} \leq 1$, in the interval $0 \leq t \leq \sqrt{\alpha - n}$
we have a subgaussian deviation inequality
\begin{equation}
\m_{2n,\alpha}\big\{(x,y)\in \R^n \times \R^n: 
\sqrt{\alpha - n}\, |f(x) - f(y)| \geq t\big\} \leq 2\,e^{-t^2/12}.\label{subgausineq}
\end{equation}
In particular, if $\alpha \geq n^2$, then 
\begin{equation}\label{subgausineq2}
\m_{n,\alpha} \otimes \m_{n,\alpha} \big\{\sqrt{\alpha - n}\, |f(x) - f(y)| \geq t\big\} \leq
4\,e^{-t^2/12}.
\end{equation}
\end{corollary}

\vskip5mm
\section{{\bf Spherical cups and associated transforms}}\label{SPHERICAL CUPS}
\setcounter{equation}{0}

\vskip2mm
\noindent
First let us employ Pisier's argument in order to prove Theorem \ref{maintheorem}.

\noindent Suppose that we are given a positive Borel measure $\mu$ on
$\R^n \times \R^n$ (in general, not necessarily finite). Given a smooth function 
$f$ on $\R^n$, we consider fluctuations of $\Delta = f(y) - f(x)$ under $\mu$. 
So, introduce the same path $x(t) = x \cos t + y \sin t$, $0 \leq t \leq \frac{\pi}{2}$, 
connecting $x$ with $y$, i.e., with $x(0) = 0$, $x(\pi/2) = y$. Again, 
$$
\Delta = \int_0^{\pi/2} \frac{d}{dt} f(x(t))\,dt =\int_0^{\pi/2}
\left<\nabla f(x(t)),x'(t)\right>dt.
$$
Hence, given a convex non-negative function $\Psi$, by Jensen's inequality,
$$
\Psi(\Delta) \leq \frac{2}{\pi} \int_0^{\pi/2} \Psi\Big(\frac{\pi}{2}
\left<\nabla f(x(t)),x'(t)\right>\Big)\,dt,
$$
and after integration, we get
\begin{equation}
\int_{\R^n} \int_{\R^n} \Psi(\Delta)\,d\mu \leq\label{firstintegpsi}
\frac{2}{\pi} \int_0^{\pi/2} \int_{\R^{2n}} \Psi\Big(\frac{\pi}{2}
\left<\nabla f(x(t)),x'(t)\right>\Big)\,dt\,d\mu.
\end{equation}
As in Definition \ref{sphericalcup}, denote by $\mu_t$ the image
of $\mu$ under the orthogonal linear transformation
$U_t: \R^n \times \R^n \rightarrow \R^n \times \R^n$ given by
$$
U_t(x,y) = (u,v) = (x(t),x'(t)) = (x \cos t + y \sin t,-x \sin t + y \cos t).
$$
Then using the change of variables formula, \eqref{firstintegpsi} may be rewritten as
$$
\int_{\R^n} \int_{\R^n} \Psi(\Delta)\,d\mu \leq
\frac{2}{\pi} \int_0^{\pi/2} \int_{\R^{2n}} \Psi\Big(\frac{\pi}{2}
\left<\nabla f(u),v\right>\Big)\,dt\,d\mu_t(u,v),
$$
which is the desired relation \eqref{maininequal} for the spherical cup $\hat \mu$ defined in \eqref{sphericalcup}.
\qed

\vskip2mm
Let us state \eqref{maininequal} once more on a functional level, assuming that $\mu$ has density
$w(x,y)$ with respect to the Lebesgue measure. From the definition of $\mu_t$
it follows that, for any bounded measurable function $g$ on $\R^n \times \R^n$,
$$
\int_{\R^n} \int_{\R^n} g(u,v)\,d\mu_t(u,v) =
\int_{\R^n} \int_{\R^n} g(u,v)\,w(u\cos t - v\sin t,u \sin t + v \cos t)\,du\, dv.
$$
Hence, according to Definition \ref{sphericalcup}, the spherical cup $\hat \mu$ has density
\begin{equation}
(U w)(u,v) = \hat w(u,v) =
\frac{2}{\pi} \int_0^{\frac{\pi}{2}} w(u\cos t - v\sin t,u \sin t + v \cos t)\,dt.\label{ptilde}
\end{equation}
As a result, we obtain the transform $U(w) = \hat w$, which may be extended 
by this formula as a linear positive operator on $L^1(\R^n \times \R^n)$.
Theorem \ref{maintheorem} has therefore the following functional counterpart.

\vskip5mm
\begin{theorem}\label{mainfunctionalcount}
Given a smooth function $f$ on $\R^n$,
for any non-negative Borel measurable function $w(x,y)$ on $\R^n \times \R^n$,
$$
\int_{\R^n} \int_{\R^n} \Psi(f(y) - f(x))\,w(x,y)\,dx\,dy \leq
\int_{\R^n} \int_{\R^n} \Psi\Big(\frac{\pi}{2} \left<\nabla f(u),v\right>\Big)\,
\hat w(u,v)\,du\, dv.
$$
In particular, for any $p \geq 1$,
$$
\int_{\R^n} \int_{\R^n} |f(y) - f(x)|^p\,w(x,y)\,dx\,dy \leq
\Big(\frac{\pi}{2}\Big)^p \int_{\R^n} \int_{\R^n} 
|\left<\nabla f(u),v\right>|^p\, \hat w(u,v)\,du\, dv.
$$
\end{theorem}

\vskip2mm
Let us stress how the operator $U$ defined in \eqref{ptilde} is acting on the 
Lebesgue spaces $L^p = L^p(\R^n \times \R^n,dx\,dy)$ with standard norms
$\|\cdot\|_p$.

\vskip5mm
{\bf Proposition 3.2.} {\sl The linear operator $U$ defined in  \eqref{ptilde} is acting 
as a contraction from $L^p$ to $L^p$ for any $p \in [1,\infty]$.
Moreover, it represents a unitary operator on $L^2$.
}

\vskip5mm
{\bf Proof.} For each $t$, the linear mapping
\begin{equation}
T_t(u,v) = (u\cos t - v\sin t,u \sin t + v \cos t), \quad u,v \in \R^n,\label{rotation}
\end{equation}
pushes forward the Lebesgue measure $\lambda_{2n}$ on $\R^n \times \R^n$
to $\lambda_{2n}$. Hence, the first claim immediately follows from \eqref{rotation}
and the triangle inequality:
$$
\|U w\|_p \leq \frac{2}{\pi} \int_0^{\frac{\pi}{2}} \|w(T_t)\|_p\,dt
 = \frac{2}{\pi} \int_0^{\frac{\pi}{2}} \|w\|_p\,dt = \|w\|_p.
$$

For the second one, let us recall that the equation $T_t(u,v) = (x,y)$ is solved as
\begin{eqnarray}
(u,v)
 & = &
T_t^{-1}(x,y) \, = \, U_t(x,y) \nonumber \\
 & = &
(x\cos t + y\sin t,-x \sin t + y \cos t) = R(T_t(y,x)),
\quad x,y \in \R^n\nonumber,
\end{eqnarray}
where $R(\xi,\eta) = (\eta,\xi)$, $\xi,\eta \in \R^n$,
is the reflection around the main diagonal in the space
$\R^n \times \R^n$. Hence, for all complex-valued functions $f,g \in L^2$,
\bee
\left<Uf,g\right>_{L^2}
 & = &
\frac{2}{\pi} \int_0^{\frac{\pi}{2}}
\bigg[\int f(T_t(u,v))\, \bar g(u,v)\,du\,dv\bigg]\,dt \\
 & = &
\frac{2}{\pi} \int_0^{\frac{\pi}{2}}
\bigg[\int f(x,y)\, \bar g(T_t^{-1}(x,y))\,dx\,dy\bigg]\,dt \\
 & = &
\frac{2}{\pi} \int_0^{\frac{\pi}{2}}
\bigg[\int f(x,y)\, \bar g(R(T_t(y,x)))\,dx\,dy\bigg]\,dt \\
 & = &
\frac{2}{\pi} \int_0^{\frac{\pi}{2}}
\bigg[\int f(R(x,y))\, \bar g(T_t(y,x))\,dx\,dy\bigg]\,dt \\
 & = &
\frac{2}{\pi} \int_0^{\frac{\pi}{2}}
\bigg[\int f(y,x)\, \bar g(T_t(y,x))\,dx\,dy\bigg]\,dt \, = \,
\left<f,Ug\right>_{L^2}. \\
\ene
The latter means that the operator $U$ is unitary.
\qed

\vskip5mm
\section{{\bf Isotropic measures}}\label{ISOTROPIC MEASURES}
\setcounter{equation}{0}

\vskip2mm
\noindent
One important particular case in Theorem \ref{maintheorem} is the bound
\begin{equation}\label{maintheoremp=2}
\int_{\R^n} \int_{\R^n} |f(y) - f(x)|^2\,d\mu(x,y) \leq \frac{\pi^2}{4}
\int_{\R^n} \int_{\R^n} \left<\nabla f(u),v\right>^2 d\hat \mu(u,v).
\end{equation}
In order to simplify the right integral, some additional assumptions
about the measures $\mu$ or $\hat \mu$ are needed. A good possible property is
isotropy. A finite Borel measure $\lambda$ on $\R^n$ is called isotropic, if
$$
\int \left<\theta,v\right>^2\,d\lambda(v) = \sigma^2 |\theta|^2 \quad
{\rm for \ all} \ \ \theta \in \R^n
$$
with some finite $\sigma^2$ ($\sigma \geq 0$) independent of $\theta$, called
the isotropic constant of $\lambda$. In this case,
\begin{equation}
\sigma^2 = \frac{1}{n} \int |v|^2\,d\lambda(v).\label{isotropconst}
\end{equation}
This definition is consistent with the one in Convex Geometry, when $\lambda$
represents the Lebesgue measure restricted to a symmetric convex body. 
The body is then called isotropic.

Any rotationally invariant measure on $\R^n$ is isotropic.
In particular, when $\lambda$ represents the normalized Lebesgue measure on 
the sphere $rS^{n-1} \subset \R^n$ of radius $r>0$ with center at the origin,
then $\sigma^2 = r^2/n$ according to \eqref{isotropconst}.
Intuitively, the isotropy means that in some sense the measure $\lambda$ is
``round'' and is not dilated in any direction. On the other hand, this is
a normalization type condition. By simple algebra, any measure with finite
second moment becomes isotropic after some linear transformation.
Note that an orthogonal transformation of an isotropic measure is isotropic.

When $\lambda$ is a probability measure on $\R^n$, and a random vector
$X = (X_1,\dots,X_n)$ is distributed according to $\lambda$, the property
of being isotropic is equivalent to the non-correlatedness of coordinates and
the requirement that all coordinates have equal $L^2$-norms:
$\E X_i X_j = \sigma^2 \delta_{ij}$.

Now, given a finite Borel measure $\nu = \nu(du,dv)$ on $\R^n \times \R^n$, let
$\pi$ denote its projection to the ``first'' coordinate $u$. In other words,
it represents the first marginal of $\nu$ defined by
$$
\pi(A) = \nu(A \times \R^n), \quad A \subset \R^n \ ({\rm Borel}).
$$
According to the general Measure Theory (cf. \cite{Roh}, \cite{Bog}), for the partition
$$
\R^n \times \R^n = \bigcup_{u \in \R^n} \{u\} \times \R^n,
$$
there exists a (unique) family of finite measures $\nu_u$ defined for
$\pi$-almost all $u$, called conditional measures, such that each $\nu_u$ is
supported on $\{u\} \times \R^n$ and
\begin{equation}
\nu = \int_{\R^n} \nu_u\,d\pi(u).\label{condmeas}
\end{equation}
Moreover, if $\nu$ is a probability measure, then $\pi$ and $\nu_u$ are
probability measures as well.

It will be more convenient to consider $\R^n$ as the space for the support
of $\nu_u$ rather then the slice $\{u\} \times \R^n$, so that
\begin{equation}\label{decompcondmeas}
\int_{\R^n} \int_{\R^n} g(u,v)\,d\nu(u,v) = 
\int_{\R^n} \bigg[\int_{\R^n} g(u,v)\, d\nu_u(v)\bigg]\,d\pi(u)
\end{equation}
for any $\nu$-integrable function $g$ on $\R^n \times \R^n$. Equality \eqref{condmeas}
thus provides a canonical representation for $\nu$ in the case of the partition of
the $2n$-space along the first coordinate, and we use ( for computations.
With this convention, if $\nu$ has density $q(u,v)$, then the marginal measure
$\pi$ has density $u \rightarrow c_u = \int q(u,v)\,dv$, and $\nu_u$ has density
$$
q_u(v) = \frac{1}{c_u}\,q(u,v) \quad {\rm for} \ \ c_u > 0.
$$

\vskip2mm
\begin{definition}\label{isotropicalongfirstcoor} Let us say that a finite measure $\nu$ on
$\R^n \times \R^n$ is isotropic along the first coordinate, if $\pi$-almost
all conditional measures $\nu_u$ in the canonical representation \eqref{condmeas} are
isotropic on $\R^n$. Equivalently, for $\pi$-almost all $u$,
$$
\int_{\R^n} \left<\theta,v\right>^2\,d\nu_u(v) = \sigma^2(u) |\theta|^2 \quad
(\theta \in \R^n)
$$
with some finite $\sigma^2(u)$, which we call the isotropic function of $\nu$
along the first coordinate.
\end{definition}
\vskip5mm
According to \eqref{isotropconst}, in that case the isotropic function is given by
\begin{equation}
\sigma^2(u) = \frac{1}{n} \int_{\R^n} |v|^2\,d\nu_u(v).\label{isotropconstfirstcoord}
\end{equation}

One can now return to Theorem 2.2, assuming that $\nu = \hat \mu$ is isotropic
along the first coordinate, with the isotropic function $\hat \sigma^2(u)$. 
Using \eqref{decompcondmeas}, on the right-hand side of \eqref{maintheoremp=2} we then deal with the integral
\bee
\int_{\R^n \times \R^n} \left<\nabla f(u),v\right>^2\, d\hat \mu(u,v)
 & = &
\int_{\R^n} \bigg[\int_{\R^n} \left<\nabla f(u),v\right>^2\, d\hat \mu_u(v)\bigg]\,d\pi(u) \\
 & = &
\int_{\R^n} \hat \sigma^2(u)\, |\nabla f(u)|^2\,d\pi(u).
\ene
As a result, we arrive in \eqref{maintheoremp=2} at the weighted Poincar\'e-type inequality.

A similar conclusion can also be made in the case $\Psi(r) = |r|$ in Theorem 2.2, 
which leads to the weighted Cheeger-type inequality. Assume that $\mu$ is
a probability measure, so that $\nu = \hat \mu$ and $\nu_u$ are probability measures as well.
Then, using \eqref{decompcondmeas} and applying Cauchy's inequality, we get
\bee
\int_{\R^n \times \R^n} \Psi\Big(\frac{\pi}{2} \left<\nabla f(u),v\right>\Big)\, d\hat \mu(u,v)
 & = &
\frac{\pi}{2} \int_{\R^n} \bigg[\int_{\R^n} |\left<\nabla f(u),v\right>|\, d\hat \mu_u(v)\bigg]\,d\pi(u) \\
 & \leq &
\frac{\pi}{2} \int_{\R^n}
\bigg(\int_{\R^n} \left<\nabla f(u),v\right>^2 d\hat \mu_u(v)\bigg)^{1/2}\,d\pi(u) \\
 & = &
\frac{\pi}{2} \int_{\R^n} \hat \sigma(u)\, |\nabla f(u)|\,d\pi(u).
\ene
Let us summarize.

\vskip5mm
\begin{corollary}\label{poincpartialisotrmeas}
Suppose that, for a finite measure $\mu$ on $\R^n \times \R^n$, 
the spherical cup $\tilde \mu$ is an isotropic measure on $\R^n \times \R^n$ along 
the first coordinate with the isotropic function $\tilde \sigma^2$. Let $\pi$ be 
the first marginal of $\hat \mu$. Then, for any smooth function $f$ on $\R^n$,
$$
\int_{\R^n \times \R^n} |f(y) - f(x)|^2\,d\mu(x,y) \leq
\frac{\pi^2}{4} \int_{\R^n} \hat \sigma^2(u)\, |\nabla f(u)|^2\,d\pi(u).
$$
In addition, if $\mu$ is a probability measure, then
$$
\int_{\R^n \times \R^n} |f(y) - f(x)|\,d\mu(x,y) \leq
\frac{\pi}{2} \int_{\R^n} \hat \sigma(u)\, |\nabla f(u)|\,d\pi(u).
$$
\end{corollary}

\vskip5mm
\section{{\bf Spherically invariant measures}}\label{SPHERICALLY INVARIANT MEASURES}
\setcounter{equation}{0}

\noindent
A natural generalization of Pisier's result is Theorem \ref{maintheorem} with measures $\mu$ 
that are spherically (rotationally) invariant on $\R^n \times \R^n$.
In the absolutely continuous case, this means that $\mu$ has a density
$$
w(x,y) = w(r) = w\big(\sqrt{|x|^2 + |y|^2}\big),
$$
depending only on the norm $r = |(x,y)|$. Since the transforms $U_t$ and 
$T_t$ are orthogonal, we have $|U_t(x,y)| = |T_t(x,y)| = |(x,y)|$. So, 
$\mu_t = \mu$ for all $t$ and thus $\hat \mu = \mu$ similarly 
to the Gaussian measures $\mu = \gamma_n \otimes \gamma_n$. 
Let us emphasize this particular case in \eqref{maininequal} once more.

\vskip5mm
\begin{theorem}\label{poincrotatinv}
Let a positive Borel measure $\mu$ on $\R^n \times \R^n$ 
be rotationally invariant. Then, for any smooth function $f$ on $\R^n$,
$$
\int_{\R^n \times \R^n} \Psi(f(y) - f(x))\,d\mu(x,y) \leq
\int_{\R^n \times \R^n} \Psi\Big(\frac{\pi}{2} \left<\nabla f(u),v\right>\Big)\,d\mu(u,v).
$$
In particular,
\begin{equation}
\int_{\R^n \times \R^n} |f(y) - f(x)|^2\,d\mu(x,y) \leq \frac{\pi^2}{4}
\int_{\R^n \times \R^n} \left<\nabla f(u),v\right>^2\,d\mu(u,v).\label{rotatinvariant}
\end{equation}
\end{theorem}

\vskip2mm
Obviously, all rotationally invariant measures are isotropic along the first
coordinate, so Corollary \ref{poincpartialisotrmeas} is applicable, and one can evaluate explicitly
the involved marginal $\pi$ and the isotropic function $\sigma^2$ which
serves as a weight (since $\hat \mu = \mu$, we may omit the hat-sign). 

Let us illustrate these relations on the example of the uniform distribution 
on the sphere. In the sequel, we denote by 
\begin{equation}
\omega_n = \frac{\pi^{\frac{n}{2}}}{\frac{n}{2}\, \Gamma(\frac{n}{2})}\label{volumeball}
\end{equation}
the volume of the unit ball $B_n = \{x \in \R^n: |x|^2 \leq 1\}$, and then 
$s_{n-1} = n\omega_n$ describes the surface measure of the unit sphere 
$S^{n-1}$. Thus, let $\mu = \sigma_{2n-1}$, $n \geq 2$,
be the normalized Lebesgue measure on the $(2n-1)$-dimensional unit sphere,
which may also be defined as the normalized restriction of the Hausdorff measure 
of dimension $2n-1$ to the sphere
$$
S^{2n-1} = \{(x,y) \in \R^n \times \R^n: |x|^2 + |y|^2 = 1\}.
$$
Every section
$$
S^{2n-1}_x = \big\{y \in \R^n: (x,y) \in S^{2n-1}\big\} =
\sqrt{1- |x|^2}\ S^{n-1} \quad (x \in \R^n, \ |x| < 1)
$$
represents the sphere in $\R^{n}$ of radius $\sqrt{1-|x|^2}$ with center
at the origin. In this case, the conditional measure $\nu_x$ in \eqref{condmeas}-\eqref{decompcondmeas} for
$\nu = \sigma_{2n-1}$ represents the uniform  distribution on $S^{2n-1}_x$.
Therefore, by \eqref{isotropconstfirstcoord}, the isotropic function for $\mu$ is given by
\begin{equation}
\sigma^2(x) = \frac{1}{n}\int_{S^{2n-1}_x} |y|^2\, d\nu_{x}(y) =
\frac{1 - |x|^2}{n}, \quad x \in \R^n, \ |x| \leq 1.\label{isotropicfunctsphere}
\end{equation}

Now, let us describe the coresponding marginal distribution
$$
\pi(A) = \sigma_{2n-1}\big(S^{2n-1}\cap(A\times\R^{n})\big), \quad
A \subset \R^n \ ({\rm Borel}).
$$
It is supported on the unit ball $B_n$, where it has a density of the form $q(x) = q(|x|)$, 
i.e. $\pi$ is also rotationally invariant. To find $q$, let $\mu_\ep$, $\ep>0$, denote the
uniform distribution in the region
$$
D_\ep = \big\{(x,y) \in \R^n \times \R^n: 1 < \sqrt{|x|^2 + |y|^2} < 1+\ep\big\}.
$$
Then $\mu_\ep \rightarrow \sigma_{2n-1}$ weakly as $\ep \rightarrow 0$, and a similar
convergence holds true for the marginal distributions of these measures, $\pi_\ep$ and $\pi$.
Moreover, the density $q_\ep(x)$ of $\pi_\ep$ is convergent to $q(x)$ on the unit ball $|x|<1$.  Denoting by ${\rm mes}_n$ the Lebesgue measure on $\R^n$, for any Borel set $A\subset B_{n}$ we have by Fubini's theorem
\[
\pi_{\ep}(A)=\frac{1}{{\rm mes}_{2n}(D_{\ep})}\int_{D_{\ep}\cap(A\times\R^{n})}dx\,dy=\frac{1}{{\rm mes}_{2n}(D_{\ep})}\int_{A}dx\int_{\left\{y\in\R^{n}:(x,y)\in D_{\ep}\right\}}dy.
\]
Therefore, for all the points $|x|<1$ we have 
\bee
q_\ep(x) 
 & = &
\frac{1}{{\rm mes}_{2n}(D_\ep)}\,{\rm mes}_n\big\{y \in \R^n:1 < \sqrt{|x|^2 + |y|^2} < 1+\ep\big\} \\
 & = &
\frac{\omega_n}{\omega_{2n}\,((1+\ep)^{2n} - 1)}\,
\Big(((1+\ep)^2 - |x|^2)^{\frac{n}{2}} - (1 - |x|^2)^{\frac{n}{2}}\Big).
\ene
Let $u(\ep)$ denote the expression in the last brackets. We have $u(0) = 0$,
$u'(0) = n(1- |x|^2)^{\frac{n}{2} - 1}$. Since $(1+\ep)^{2n} - 1 = 2n \ep + O(\ep^2)$,
it follows that
$$
q_\ep(x) = \frac{\omega_n}{2 \omega_{2n}}\,(1- |x|^2)^{\frac{n}{2} - 1} + O(\ep)
$$
and therefore
$$
q(x) = \frac{\omega_n}{2 \omega_{2n}}\,(1- |x|^2)^{\frac{n}{2} - 1}, \quad |x|<1.
$$
To describe the coefficient explicitly, one may refer to the formula \eqref{volumeball} which gives
$$
\frac{n\omega_n}{2n\, \omega_{2n}} = 
\frac{\Gamma(n)}{\pi^{\frac{n}{2}}\, \Gamma(\frac{n}{2})}.
$$
Thus, the marginal distribution of $\sigma_{2n-1}$ is the probability measure 
on the unit ball $B_n$ in $\R^n$ with density
\begin{equation}
\frac{d\pi(u)}{du} = \frac{\Gamma(n)}{\pi^{\frac{n}{2}}\, \Gamma(\frac{n}{2})}\,
(1 - |u|^2)^{\frac{n}{2}-1}, \quad |u|<1.\label{marginaldistrsphere}
\end{equation}
Hence, recalling  formula \eqref{isotropicfunctsphere} for the isotropic function,
we can apply and state Corollary \ref{poincpartialisotrmeas} in a more explicit form.

\vskip5mm
\begin{corollary}\label{corpoincineqsphere}
For any smooth function $f$ on $\R^n$,
\begin{equation}\label{poincineqsphere}
\int_{S^{2n-1}} |f(x) - f(y)|^2\,d\sigma_{2n-1}(x,y) \, \leq \,
\frac{\pi^2}{4n}
\int_{|u| < 1} |\nabla f(u)|^2\, (1 - |u|^2)\,d\pi(u),
\end{equation}
where the probability measure $\pi$ is given in \eqref{marginaldistrsphere}.
\end{corollary}

\vskip5mm
Let us see what kind of concentration is hidden in \eqref{poincineqsphere}. As easy to check,
for any $p>0$,
$$
\int_{|u| < 1} (1 - |u|^2)^{p-1}\,du = 
\pi^{\frac{n}{2}}\,\frac{\Gamma(p)}{\Gamma(\frac{n}{2} + p)}.
$$
We use this identity for $p = \frac{n}{2} + 1$.
If $|\nabla f| \leq 1$, the last integral in \eqref{poincineqsphere} does not exceed
$$
\int_{|u| < 1} (1 - |u|^2)\,d\pi(u) = 
\frac{\Gamma(n)}{\pi^{\frac{n}{2}}\, \Gamma(\frac{n}{2})}\,
\frac{\Gamma(p)}{\Gamma(\frac{n}{2} + p)} = \frac{1}{2}.
$$
Hence, by \eqref{poincineqsphere},
$$
\int_{S^{2n-1}} |f(x) - f(y)|^2\,d\sigma_{2n-1}(x,y) \, \leq \,
\frac{\pi^2}{8n}.
$$
A similar conclusion can be made on the basis of the Poincar\'e
inequality on the sphere $S^{2n-1}$.

The inequality (5.5) is rather similar to the Poincar\'e-type inequality
for the measure $\pi$,
$$
\int_{S^{2n-1}} |f(x) - f(y)|^2\,d\pi(x)\,d\pi(y) \, \leq \,
\frac{C}{n} \int_{|u| < 1} |\nabla f(u)|^2\,d\pi(u),
$$
which was derived in \cite{B2} in a more general setting
of log-concave, spherically symmetric probability measures on $\R^n$.

\vskip5mm
\section{{\bf Background on Cauchy measures}}\label{BACKGROUND CAUCHY}
\setcounter{equation}{0}

\vskip2mm
\noindent
Let us recall basic definitions and facts about the multidimensional Cauchy 
distributions. 

\vskip2mm
\begin{definition}\label{Cacuhymeas} The $n$-dimensional Cauchy measure $\m_{n,\alpha}$ on $\R^n$ 
of order $\alpha>0$ has density
\begin{equation}
w_{n,\alpha}(x) = 
\frac{1}{c_{n,\alpha}}\,(1 + |x|^2)^{-\alpha}, \quad x \in \R^n,\label{cauchyweight}
\end{equation}
where $c_{n,\alpha}$ is a normalizing constant such that $\m_{n,\alpha}(\R^n) = 1$.
\end{definition}
\vskip2mm
To see when the function $w_\alpha(x)$ is integrable, one may 
integrate in polar coordinates:
\begin{equation}
\int_{\R^n} \frac{dx}{(1 + |x|^2)^\alpha} =
n \omega_n \int_0^\infty \frac{r^{n-1}}{(1 + r^2)^\alpha}\,dr.\label{integnormalizconst}
\end{equation}
The latter integral is finite if and only if $\alpha > \frac{n}{2}$. This will be always
assumed when speaking about Cauchy measures on $\R^n$. 

Changing the variable $r = \sqrt{s}$ and 
then $s = \frac{1}{t} - 1$, the last integral in \eqref{integnormalizconst} is equal to
\bee
\frac{1}{2} \int_0^\infty \frac{s^{\frac{n}{2} - 1}}{(1 + s)^\alpha}\,ds 
 & = &
\frac{1}{2} \int_0^1 t^{\alpha - \frac{n}{2} - 1}\, (1-t)^{\frac{n}{2} - 1}\,dt \\ 
& = &
\frac{1}{2} \, B\Big(\alpha - \frac{n}{2},\frac{n}{2}\Big) \, = \,
\frac{1}{2} \, \frac{\Gamma(\alpha - \frac{n}{2})\,\Gamma(\frac{n}{2})}{\Gamma(\alpha)}.
\ene
Hence, the normalizing constant in is given by
\begin{equation}\label{normalizconst}
c_{n,\alpha} = \frac{n\, \omega_n}{2} \, 
\frac{\Gamma(\alpha - \frac{n}{2})\,\Gamma(\frac{n}{2})}{\Gamma(\alpha)} = 
\pi^{\frac{n}{2}} \, \frac{\Gamma(\alpha - \frac{n}{2})}{\Gamma(\alpha)}.
\end{equation}

\vskip2mm
{\bf Characterization by means of $\chi$ squared distribution.} 
Write $\alpha = \frac{n+d}{2}$ for 
$$
d = 2\alpha - n > 0. 
$$
The probability measure $\m_{n,\alpha}$ may be characterized as 
the distribution of the random vector $X = Z/\eta$, where $Z = (Z_1,\dots,Z_n)$
is a random vector in $\R^n$ with the standard Gaussian distribution $\gamma_n$, 
and where $\eta > 0$ is a random variable independent of $Z$ and having the 
$\chi_d$ -distribution with $d$ degrees of freedom. That is, $\eta$ has density
$$
\chi_d(r) = \frac{1}{2^{\frac{d}{2} - 1}\,\Gamma(\frac{d}{2})}\,r^{d-1}\,e^{-r^2/2},
\quad r> 0.
$$

Hence, $\m_{n,\alpha}$ represents the image of the product measure
$\gamma_n \otimes \chi_d$ on $\R^n \times (0,\infty)$ under the map 
$(z, r) \rightarrow z/r$ (with some abuse, we denote by $\chi_d$ 
the measure with density $\chi_d(r)$).
When $d$ is integer, $\chi_d$ represents the distribution
of the Euclidean norm $\eta = |V| = (\xi_1^2 + \dots + \xi_d^2)^{1/2}$ of 
a standard normal random vector $V = (\xi_1,\dots,\xi_d)$ in $\R^d$.

\vskip2mm
{\bf Essential support of $\m_{n,\alpha}$.}
Note that $\E\,|Z|^2 = n$, while $\E\,|V|^2 = d$. Moreover,
for large $n$ and $d$, with high probability we have $\frac{n}{2} < |Z|^2 < 2n$
and $\frac{d}{2} < |V|^2 < 2d$ (as a consequence of the Gaussian concentration
phenomenon). Hence, with high probability
$$
\frac{n}{4d} < |X|^2 < \frac{4n}{d}.
$$
In other words, the measure $\m_{n,\alpha}$ is almost concentrated on the
Euclidean ball with center at zero and radius $R$ of order $\sqrt{n/d}$. If $\alpha$ is
essentailly greater than $\frac{n}{2}$, for example, if $\alpha \geq n$,
then $d$ is of order $2\alpha$, so that $R$ is approximately  $\sqrt{n/(2\alpha)}$.

\vskip2mm
{\bf Projections.} The $\chi_d$-characterization has a useful consequence 
about the images $\m_{n,\alpha} T^{-1}$ of $\m_{n,\alpha}$ under linear 
projections $T:\R^n \rightarrow H$ to $k$-dimensional linear subspaces $H$ 
in $\R^n$ ($1 \leq k \leq n$). Since $\m_{n,\alpha}$ is spherically invariant,
it is sufficient to consider canonical projections from $\R^n$ to $\R^k$,
that is, the maps
$$
T(x_1,\dots,x_n) = (x_1,\dots,x_k), \quad (x_1,\dots,x_n) \in \R^n.
$$

Let $Z = (Z_1,\dots,Z_n)$ be a random vector in $\R^n$ with the standard 
Gaussian distribution, and $\eta$ be a random variable independent of 
$Z$ and having the $\chi_d$ -distribution with $d = 2\alpha - n$ degrees of freedom. 
Then the random vector $X = Z/\eta$ has distribution $\m_{n,\alpha}$.
Moreover,
$$
T(X) = \frac{1}{\eta}\,T(Z) = \frac{1}{\eta}\,(Z_1,\dots,Z_k)
$$
is of a similar form as $X$ itself. Since $(Z_1,\dots,Z_k)$ is a standard normal
random vector in $\R^k$, we conclude that the projection
$\m_{n,\alpha} T^{-1}$ represents the $k$-dimensional Cauchy measure
of the  order $\beta$ such that $d = 2\beta - k$. That is, the $k$-dimensional
projection of $\m_{n,\alpha}$ is a Cauchy measure $\m_{k,\beta}$
of order
\begin{equation}
\beta = \frac{d + k}{2} = \alpha - \frac{n-k}{2}.\label{beta}
\end{equation}

\vskip2mm
{\bf Gaussian limit.} If a random vector $X$ has a Cauchy distribution 
$\m_{n,\alpha}$, then the random vector $X_\alpha = \sqrt{2\alpha}\,X$ has density
$$
\widetilde w_{n,\alpha}(x) = (2\alpha)^{-n/2}\,
w_{n,\alpha}\Big(\frac{1}{\sqrt{2\alpha}}\,x\Big) = 
\frac{1}{c_{n,\alpha}' (1 + \frac{1}{2\alpha}\, |x|^2)^\alpha}, \quad x \in \R^n,
$$
where $c_{n,\alpha'} = (2\alpha)^{n/2}\,c_{n,\alpha}$. This normalization is consistent
with the size of the essential support of the Cauchy measure for large values of
$\alpha$ (in which case the values of $|X_\alpha|$ do not exceed a muptiple of
$\sqrt{n}$ in a full analogy with standard Gaussian random vectors). Indeed,
$$
\widetilde w_{n,\alpha}(x) \rightarrow (2\pi)^{-n/2}\,e^{-|x|^2/2} \quad {\rm as} \ \ 
\alpha \rightarrow \infty.
$$
Therefore, the linear image $\widetilde \m_{n,\alpha}$ of $\m_{n,\alpha}$ under 
the map $x \rightarrow \sqrt{2\alpha}\, x$ approaches for growing $\alpha$ the 
standard Gaussian measure $\gamma_n$ on $\R^n$ in total variation norm.
This property also follows from the above $\chi_d$-characterization and the
convergence in distribution
$$
\frac{|V|^2}{d} = \frac{\xi_1^2 + \dots + \xi_d^2}{d} \rightarrow 1 \quad {\rm as} \ \ 
d \rightarrow \infty,
$$
which is the weak law of large numbers for independent 
normal random variables $\xi_i$.

Thus, the class of Cauchy measures may serve as a pre-Gaussian model,
in the sense that many properties of the Gaussian measure may be potentially
obtained in the limit from similar ones for the Cauchy measures. 

\vskip2mm
{\bf Moments of the Euclidean norm.} Since the density of the Cauchy measure 
$\m_{n,\alpha}$ on $\R^n$ decay polynomially at infinity, the Euclidean norm has finite
$L^p$-norms for a certain range $p \geq 0$ only, namely, for $p < 2\alpha - n$. 
Integrating in polar coordinates and repeating the previous computations, we have
$$
\int_{\R^n} |x|^p\,dm_\alpha(x) =
\frac{1}{c_{n,\alpha}} \int_{\R^n} \frac{|x|^p}{(1 + |x|^2)^\alpha}\,dx = 
\frac{\Gamma(\alpha - \frac{n+p}{2})\,
\Gamma(\frac{n+p}{2})}{\Gamma(\alpha - \frac{n}{2})\,\Gamma(\frac{n}{2})}.
$$
In particular, the second moment is finite if and only if
$\alpha > \frac{n+2}{2}$, in which case
$$
\int_{\R^n} |x|^2\,d\m_{n,\alpha}(x) = \frac{n}{2\alpha - n - 2}, \quad
\int_{\R^n} \left<x,\theta\right>^2 d\m_{n,\alpha}(x) = \frac{|\theta|^2}{2\alpha - n - 2}.
$$
This also shows that the $\sqrt{2\alpha}$-normalization in the definition
of $\widetilde \m_{n,\alpha}$ is natural, when $\alpha$ is essentialy larger than $n$.

\vskip5mm
\section{{\bf Cauchy measures on $\R^n \times \R^n$. Proof of Theorem \ref{Poinccauchytheo}}}\label{PROOF THEOREM CAUCHY}
\setcounter{equation}{0}

\vskip2mm
\noindent
According to the general definition \eqref{cauchyweight} and the formulas \eqref{volumeball}--\eqref{normalizconst}, the 
$2n$-dimensional Cauchy measures $\m_{2n,\alpha}$ on $\R^n \times \R^n$ 
have densities
\begin{equation}\label{densityCauchyRnRn}
w_{2n,\alpha}(x,y) = \frac{1}{c_{2n,\alpha}}\, 
(1 + |x|^2+ |y|^2)^{-\alpha}, \quad x,y \in \R^n,
\end{equation}
where $\alpha>n$ is a parameter and where the normalizing constant is given by
\begin{equation}\label{normalizconstRnRn}
c_{2n,\alpha} = n\, \omega_{2n} \, 
\frac{\Gamma(\alpha - n)\,\Gamma(n)}{\Gamma(\alpha)} = 
\pi^n \, \frac{\Gamma(\alpha - n)}{\Gamma(\alpha)}.
\end{equation}

As we know from \eqref{beta} with $k=n$ and with $n$ being replaced with $2n$,
the marginal of $\m_{2n,\alpha}$ represents the $n$-dimensional Cauchy measure
$\m_{n,\beta}$ of order $\beta = \alpha - \frac{n}{2}$. Equivalently, if the couple 
of random vectors $(X,Y)$ in $\R^{2n}$ has a Cauchy distribution of order $\alpha$, 
then both $X$ and $Y$ have a Cauchy distribution in $\R^n$ of order $\alpha - \frac{n}{2}$. 
This fact may be obtained directly by computation of the marginal density.

By Theorem \ref{poincrotatinv} applied with $\Psi(t) = |t|^p$, $p \geq 1$, for any 
smooth function $f$ on $\R^n$,
\begin{equation}\label{ineqcauchy1}
\int_{\R^n} \int_{\R^n} |f(x) - f(y)|^p\,d\m_{2n,\alpha}(x,y) \leq \Big(\frac{\pi}{2}\Big)^p
\int_{\R^n} \int_{\R^n} |\left<\nabla f(x),y\right>|^p\,d\m_{2n,\alpha}(x,y).
\end{equation}
In order to simplify the last integral, let us fix a real number $p \geq 0$, a vector 
$v = r\theta$, $r \geq 0$, $\theta \in S^{n-1}$, and consider the integrals
\begin{equation}
I_p(x,v) =
\int_{\R^n} \frac{|\left<v,y\right>|^p}{(1 + |x|^2+ |y|^2)^\alpha}\,dy =
r^p \int_{\R^n} \frac{|\left<\theta,y\right>|^p}{(1 + |x|^2+ |y|^2)^\alpha}\,dy.\label{Ip}
\end{equation}
In particular, the quantity $I_0(x) = I_0(x,v)$ corresponding to $p=0$ 
does not depend on $v$. In this case, according to \eqref{densityCauchyRnRn}, the equality
\begin{equation}\label{nu}
\frac{d\nu(x)}{dx} = \frac{1}{c_{2n,\alpha}}\, I_0(x), \quad x \in \R^n,
\end{equation}
defines a probability measure on $\R^n$ which is the image $\nu$ of $\m_{2n,\alpha}$ under the projection
$T(x_1,\dots,x_{2n}) = (x_1,\dots,x_n)$ from $\R^n \times \R^n$ to $\R^n$.
So, $\nu$ is a Cauchy distribution on $\R^n$ of order 
$\alpha - \frac{n}{2}$.

We will recover this fact as part of a more general case $p \geq 0$ in 
\eqref{ineqcauchy1}. Put
\begin{equation}\label{G}
G(n,p) = \frac{\E\,|\xi|^p}{\E\,|Z|^p},
\end{equation}
where $\xi$ and $Z$ are standard normal random vectors in $\R$ and $\R^n$,
respectively.

\vskip5mm
\begin{lemma}\label{computIp}
 For any $\alpha> \frac{n+p}{2}$, $p \geq 0$, we have
\begin{equation}\label{expressIp}
I_p(x,v) = I_p(x)\,|v|^p, \quad \text{where }
I_p(x) = A\,(1 + |x|^2)^{-\beta}, \quad \beta = \alpha - \frac{n+p}{2},
\end{equation}
and\begin{equation}
A = G(n,p)\,
\frac{n \omega_n}{2} \, B\Big(\alpha - \frac{n+p}{2},\frac{n+p}{2}\Big).\label{constantIp}
\end{equation}
\end{lemma}

In particular,
$$
I_0(x) = \frac{A}{(1 + |x|^2)^{\alpha - \frac{n}{2}}}.
$$
This shows that the measure $\nu$ defined in \eqref{nu} is a Cauchy distribution
of order $\alpha - \frac{n}{2}$.

\vskip5mm
{\bf Proof.} By the rotational invariance of the Lebesgue 
measure, the last integral in \eqref{ineqcauchy1} does not depend on $\theta$, so that
$$
I_p(x) \equiv \E_\theta I_p(x,\theta)
$$
for all $x \in \R^n$, where $\E_\theta$ means the integral (the mean) 
with respect to $\sigma_{n-1}$. Thus, for all $\theta$,
$$
I_p(x,\theta) = 
\int_{\R^n} \frac{\E_\theta |\left<\theta,y\right>|^p}{(1 + |x|^2+ |y|^2)^\alpha}\,dy.
$$
This expectation is actually a multiple of $|y|^p$ by the rotational invariance of
$\sigma_{n-1}$, that is,
$$
\E_\theta |\left<\theta,y\right>|^p = G(n,p)\,|y|^p
$$
with some constant $G(n,p)$.
Taking $y = (1,0,\dots,0)$, we get $G(n,p) = \E_\theta |\theta_1|^p$.
Moreover, let $Z = (Z_1,\dots,Z_n)$ be a standard normal vector in $\R^n$, so that
$Z/|Z|$ is uniformly distributed in the unite sphere $S^{n-1}$ and is independent
of $|Z|$. Hence, we have the product of two independent random variables
$Z_1 = \frac{Z_1}{|Z|}\,|Z|$, which gives
$$
\E\,|Z_1|^p = \E\,\Big|\frac{Z_1}{|Z|}\Big|^p \, \E\,|Z|^p =
\E_\theta |\theta_1|^p \, \E\,|Z|^p = G(n,p) \, \E\,|Z|^p.
$$
Thus, the constant $G(n,p)$ is described according to \eqref{G}.
As a consequence,
\bee
I_p(x,v) 
 & = &
G(n,p)\, |v|^p \int_{\R^n} \frac{|y|^p}{(1 + |x|^2+ |y|^2)^\alpha}\,dy \\
 & = &
\frac{G(n,p)}{(1 + |x|^2)^{\alpha - \frac{n+p}{2}}}\, |v|^p 
\int_{\R^n} \frac{|z|^p}{(1 + |z|^2)^\alpha}\,dz.
\ene
The last integral is finite if and only if $\alpha > \frac{n+p}{2}$, which is assumed.
In this case, the last integral can be evaluated in polar coordinates: It is equal to
\bee
n \omega_n \int_0^\infty \frac{r^{n+p-1}}{(1 + r^2)^\alpha}\,dr
 & = &
\frac{n \omega_n}{2} \int_0^\infty \frac{s^{\frac{n+p}{2} - 1}}{(1 + s)^\alpha}\,ds \\
 & = &
\frac{n \omega_n}{2} \int_0^\infty t^{\alpha - \frac{n+p}{2} - 1}\, (1-t)^{\frac{n+p}{2} - 1}\,dt \\ 
& = &
\frac{n \omega_n}{2} \, B\Big(\alpha - \frac{n+p}{2},\frac{n+p}{2}\Big),
\ene
where we changed the variable $r = \sqrt{s}$ and then $s = \frac{1}{t} - 1$,
and where $B$ denotes the classical beta-function. As a result, we arrive at \eqref{expressIp} with the
constant $A$ as in \eqref{constantIp}.
\qed

\vskip5mm
{\bf Proof of Theorem \ref{Poinccauchytheo}} Using \eqref{densityCauchyRnRn} and applying Lemma \ref{computIp}, the second
integral in \eqref{ineqcauchy1} is simplied to
\bee
\frac{1}{c_{2n,\alpha}} \int_{\R^n} I_p(x,\nabla f(x))\,dx
 & = &
\frac{A}{c_{2n,\alpha}} \int_{\R^n} |\nabla f(x))|^{p}(1+|x|^{2})^{-\beta}\,dx \\
 & = &
A\,\frac{c_{n,\beta}}{c_{2n,\alpha}} \int_{\R^n} |\nabla f(x))|^{p}\,d\m_{n,\beta}.
\ene

Note that the probability measure $\m_{n,\beta}$
is defined for $\beta = \alpha - \frac{n+p}{2} > \frac{n}{2}$, which forces
$p < 2 (\alpha - n)$ given by assumption in the statement. It remains to simplify the constant in front of the last integral. Recall that,
according to \eqref{normalizconst} and \eqref{normalizconstRnRn},
\bee
c_{n,\beta}
 & = &
\frac{n \omega_n}{2} \, 
\frac{\Gamma(\beta - \frac{n}{2})\,\Gamma(\frac{n}{2})}{\Gamma(\beta)} \, = \, 
\pi^{\frac{n}{2}} \, 
\frac{\Gamma(\alpha - \frac{2n+p}{2})}{\Gamma(\alpha - \frac{n+p}{2})}, \\
c_{2n,\alpha} 
 & = &
n\, \omega_{2n} \, \frac{\Gamma(\alpha - n)\,\Gamma(n)}{\Gamma(\alpha)} \, = \,
\pi^n \, \frac{\Gamma(\alpha - n)}{\Gamma(\alpha)}.
\ene
Hence
\begin{equation}\label{ratio}
\frac{c_{n,\beta}}{c_{2n,\alpha}} = \frac{1}{\pi^{\frac{n}{2}}} \, 
\frac{\Gamma(\alpha - \frac{2n+p}{2})}{\Gamma(\alpha - \frac{n+p}{2})}\,
\frac{\Gamma(\alpha)}{\Gamma(\alpha - n)}.
\end{equation}

Next, in order to compute the constant $G(n,p)$ in \eqref{G}, we use the formula 
\eqref{volumeball} for $\omega_n$ in terms of the Gamma function and integrate in polar coordinates. 
Changing the variable $r = \sqrt{2s}$, we get that
\bee
\E\,|Z|^p
 & = &
\frac{1}{(2\pi)^{n/2}} \int_{\R^n} |x|^p\,e^{-|x|^2/2}\,dx \, = \, 
\frac{n \omega_n}{(2\pi)^{n/2}} \int_0^\infty r^{n+p-1}\,e^{-r^2/2}\,dr \\
 & = &
\frac{n \omega_n}{\pi^{n/2}}\ 2^{p/2} \int_0^\infty 
s^{\frac{n+p-1}{2}} e^{-s}\,ds \, = \, 2^{\frac{p}{2}}\,
\frac{\Gamma(\frac{n+p}{2})}{\Gamma(\frac{n}{2})}.
\ene
In particular, in dimension one 
$$
\E\,|\xi|^p \, = \, \frac{2^{\frac{p}{2}}}{\sqrt{\pi}}\,
\Gamma\Big(\frac{p+1}{2}\Big).
$$
It follows that
$$
G(n,p) = \frac{\E\,|\xi|^p}{\E\,|Z|^p} = \frac{1}{\sqrt{\pi}}\,
\frac{\Gamma(\frac{n}{2})\Gamma(\frac{p+1}{2})}{\Gamma(\frac{n+p}{2})}.
$$
Hence, according to \eqref{constantIp},
$$
A = G(n,p)\,
\frac{n \omega_n}{2} \, B\Big(\alpha - \frac{n+p}{2},\frac{n+p}{2}\Big) =
\frac{\pi^{\frac{n}{2}}}{\sqrt{\pi}}\,
\frac{\Gamma(\alpha - \frac{n+p}{2})\Gamma(\frac{p+1}{2})}{\Gamma(\alpha)}.
$$
Thus, using \eqref{expressIp},
$$
A \, \frac{c_{n,\beta}}{c_{2n,\alpha}} = 
\frac{1}{\sqrt{\pi}}\,
\frac{\Gamma(\frac{p+1}{2})\Gamma(\alpha - \frac{2n + p}{2})}{\Gamma(\alpha - n)}.
$$
\qed

\vskip5mm
\section{{\bf Poincar\'e-type inequalities for $L^1$-norms and isoperimetry}}\label{POINCARE AND ISOPERIMETRY}
\setcounter{equation}{0}

\vskip2mm
\noindent
In the important particular case $p = 1$, Theorem \eqref{Poinccauchytheo} is reduced to the following
assertion about Cauchy measures with parameters $\alpha$ on $\R^{2n}$ and 
$\beta = \alpha - \frac{n+1}{2}$ on $\R^n$.

\vskip5mm
\begin{corollary}\label{PoincareCauchyp=1}
Let $\alpha > n + \frac{1}{2}$. For any smooth function $f$ on $\R^n$,
\begin{equation}
\label{PoincareCauchyineqp=1}
\int_{\R^n \times \R^n} |f(x) - f(y)|\,d\m_{2n,\alpha}(x,y) \leq
\frac{\sqrt{\pi}}{2}\, \frac{\Gamma(\alpha - n - \frac{1}{2})}{\Gamma(\alpha - n)}
\int_{\R^n} |\nabla f| \, d\m_{n,\beta}.
\end{equation}
In particular, for $\alpha \geq n+1$
\begin{equation}
\int_{\R^n \times \R^n} |f(x) - f(y)|\,d\m_{2n,\alpha}(x,y) \leq
\sqrt{\pi}\, \frac{1}{\sqrt{\alpha - n}} \int_{\R^n} |\nabla f| \, d\m_{n,\beta}.\label{PoincareCauchyineqp=1partform}
\end{equation}
\end{corollary}

\vskip2mm
To bound from above the constant in \eqref{PoincareCauchyineqp=1} by a simpler expression, 
one may use Wendel's inequality
$\Gamma(x + \frac{1}{2}) \leq \Gamma(x) \sqrt{x}$ ($x>0$), which for $x \geq 1$ gives
$$
\Gamma\Big(x - \frac{1}{2}\Big) = \frac{1}{x - \frac{1}{2}}\,
\Gamma\Big(x + \frac{1}{2}\Big) \leq 
\frac{\sqrt{x}}{x - \frac{1}{2}}\,\Gamma(x) \leq \frac{2}{\sqrt{x}}\,\Gamma(x).
$$
Applying this with $x = \alpha - n$, we get \eqref{PoincareCauchyineqp=1partform}.

If $\alpha$ is sufficiently large, for example, $\alpha \geq 2n$, 
the constants in these inequalities do not exceed a multiple of $1/\sqrt{\alpha}$. 
Moreover, as was already explained in Section \ref{EXTENSIONS GENERAL MEASURES}, after rescaling of the space variable
and in the limit as $\alpha \rightarrow \infty$, \eqref{PoincareCauchyineqp=1} yields 
the $L^1$-Poincare-type inequality over the Gaussian measure
\begin{equation}\label{GausPoincp=1}
\int_{\R^n} \int_{\R^n} |f(x) - f(y)|\,d\gamma_n(x)\, d\gamma_n(y) \, \leq \,
\sqrt{\frac{\pi}{2}} \int_{\R^n} |\nabla f| \, d\gamma_n.
\end{equation}
Here the constant is optimal and is attained asymptotically, when $f$ approaches 
the indicator function of a half-space whose boundary passes through the origin.

Let us comment on the geometric meaning of Sobolev-type inequalities such as 
\eqref{PoincareCauchyineqp=1}-\eqref{GausPoincp=1}. One can prove the following general isoperimetric-type characterization,
using the notion of the $\nu$-perimeter $\nu^+(A)$ defined in \eqref{perimeter}.
Denote by $ A^{c} = \R^n \setminus A$ the complement of the set $A$.

\vskip5mm
\begin{lemma}\label{geometricmeaningSobineq} Let $\nu$ be an absolutely continuous probability 
measure on $\R^n$ and let $\mu$ be a finite measure on $\R^n \times \R^n$ 
which is invariant under the map $(x,y) \rightarrow (y,x)$. 
The following two assertions are equivalent:
$a)$ For any smooth function $f$ on $\R^n$,
\begin{equation}
\int_{\R^n \times \R^n} |f(x) - f(y)|\,d\mu(x,y) \leq \int_{\R^n} |\nabla f|\,d\nu.\label{Poincaremunu}
\end{equation}

$b)$ For any closed set $A$ in $\R^n$,
\begin{equation}
\nu^+(A) \geq 2\mu(A \times A^{c}).\label{isoperimetricmunu}
\end{equation}
\end{lemma}

{\bf Proof.} The argument is standard, and we give it here for completeness.
The property $a)$ can be equivalently stated for different classes of functions:

$a')$ The relation \eqref{Poincaremunu} holds for all $C^\infty$-smooth functions $f$ on $\R^n$
with a compact support;

$a'')$ The relation\eqref{Poincaremunu} holds for all locally Lipschitz functions $f$ on $\R^n$
with the generalized modulus of gradient
$$
|\nabla f(x)| = \limsup_{y \rightarrow x}\,\frac{|f(x) - f(y)|}{|x-y|}, \quad
x \in \R^n.
$$

The property of being locally Lipschitz means that $f$ has a finite
Lipschitz semi-norm in some neighborhood of any point in $\R^n$. In particular,
$f$ has to be continuous and a.e. differentiable (by Rademacher's theorem).
Moreover, $|\nabla f(x)|$ is finite, Borel measurable, and coincides 
with the modulus of the usual gradient at every point $x$ where $f$ is differentiable.

It shoud be clear that  $a'') \Rightarrow a) \Rightarrow a')$. Using a smoothing
argument, these implications can be reversed; cf. \cite{B-C-G}, Proposition 5.4.1,
for the proof of a similar statement about Poincar\'e-type inequalities for $L^2$-norms
(at this point, the absolute continuity of $\nu$ is required).

For the implication $a'') \Rightarrow b)$, given a non-empty closed set $A$ in 
$\R^n$, one may pick up functions $f_k : \R^n \rightarrow [0,1]$ with finite 
Lipschitz semi-norm such that $f_k \rightarrow 1_A$ pointwise and
$$
\limsup_{k \rightarrow \infty} \int_{\R^n} |\nabla f_k|\,d\mu \leq \nu^+(A)
$$
(cf. Proposition 5.2.2 in \cite{B-C-G}). Since, by the Lebesgue dominated convergence theorem,
\bee
\lim_{k \rightarrow \infty}
\int_{\R^n \times \R^n} |f_{k}(x) - f_{k}(y)|\,d\mu(x,y) 
 & = &
\int_{\R^n \times \R^n} |1_A(x) - 1_A(y)|\,d\mu(x,y) \\
 & = &
\mu(A \times A^{c}) + \mu(A^{c} \times A) = 2\mu(A \times  A^{c}),
\ene
the desired inequality \eqref{isoperimetricmunu} follows.

Finally, to derive \eqref{Poincaremunu} from \eqref{isoperimetricmunu}, assume that $f$ is smooth and bounded.
We may assume that $f \geq 0$, since the inequality \eqref{Poincaremunu} does not change
when adding any constant to $f$. The sets $A_t = \{x \in \R^n: f(x) \geq t\}$, $t>0$, 
are closed, so that \eqref{isoperimetricmunu} is applicable. Applying the co-area inequality 
(cf. Proposition 5.2.2 in \cite{B-C-G}) together with Fubini's theorem, we then get
\bee
\int_{\R^n} |\nabla f|\,d\nu 
 & \geq &
\int_0^\infty \nu^+(A_t)\,dt  \, \geq \,
2 \int_0^\infty \mu(A_t \times A_{t}^{c})\,dt \\
 & = &
\int_0^\infty \bigg[\int_{\R^n \times \R^n} |1_{A_t}(x) - 1_{A_t}(y)|\,d\mu(x,y)\bigg]\,dt \\
 & = &
\int_{\R^n \times \R^n} \bigg[\int_0^\infty |1_{A_t}(x) - 1_{A_t}(y)|\,dt\bigg] d\mu(x,y).
\ene
Using the triangle inequality, we see that
the last expression is greater than or equal to
$$
\int_{\R^n \times \R^n} \Big|\int_0^\infty (1_{A_t}(x) - 1_{A_t}(y))\,dt\Big|\,d\mu(x,y) = 
\int_{\R^n \times \R^n} |f(x) - f(y)|\,d\mu(x,y).
$$
As a result, we arrive at \eqref{Poincaremunu}.
\qed

\vskip2mm
When $\mu$ is a multiple of the product measure, Lemma \ref{geometricmeaningSobineq} is reduced
to the next well-known characterization.

\vskip5mm
\begin{lemma}\label{geometricmeaningSobineqproductmeas}
Let $\nu$ be an absolutely continuous probability 
measure on $\R^n$. Given $h>0$, the following two assertions are equivalent:

$a)$ For any smooth function $f$ on $\R^n$,
\begin{equation}
h \int_{\R^n} \int_{ \R^n} |f(x) - f(y)|\,d\nu(x)\,d\nu(y) \leq \int_{\R^n} |\nabla f|\,d\nu.\label{Poincaremunuproductmeas}
\end{equation}

$b)$ For any closed set $A$ in $\R^n$,
\begin{equation}
\nu^+(A) \geq 2h\,\nu(A)(1 - \nu(A)).\label{isoperimetricmunuproductmeas}
\end{equation}
\end{lemma}

As an equivalent functional form for the isoperimetric inequality \eqref{isoperimetricmunuproductmeas},
one may also consider the Sobolev-type inequality
\begin{equation}\label{Poincineqnumean}
h \int_{ \R^n} |f - m|\,d\nu \leq \int_{ \R^n} |\nabla f|\,d\nu,
\end{equation}
where $m = \int f\,d\nu$ is the $\nu$-mean of $f$, and where 
the indicator functions $f = 1_A$ still play an extremal role in the asymptotic sense. 

Both \eqref{Poincaremunuproductmeas} and \eqref{Poincineqnumean} are particular cases of Sobolev-type inequalities of the form
\begin{equation}
L f \leq \int_{ \R^n} |\nabla f|\,d\nu, \quad Lf = \sup_{g \in G} \int_{ \R^n} fg\,d\nu,\label{supfunctional}
\end{equation}
for arbitrary families $G$ of functions $g$ on $\R^n$ such that the functional
$Lf$ is well-defined (at least, for bounded $f$). Indeed, \eqref{Poincineqnumean} corresponds 
to \eqref{supfunctional} for the class $G$ of all Borel measurable $g$ with $\nu$-mean 
zero such that $|g(x)| \leq h$ for all $x \in \R^n$. Similarly, \eqref{Poincaremunuproductmeas} corresponds 
to \eqref{supfunctional} for the class $G$ of  functions $g$ representable as
$$
g(x) = \int (u(x,y) -  u(y,x))\,d\nu(y)
$$
with arbitrary Borel measurable functions $u$ such that $|u(x,y)| \leq h$, 
$x,y \in \R^n$. In this case,
\bee
\int fg\,d\nu 
 & = &
\int f(x) \bigg[\int u(x,y)\,d\nu(y)\bigg]\,d\nu(x) -
\int f(x) \bigg[\int u(y,x)\,d\nu(y)\bigg]\,d\nu(x) \\
 & = &
\int f(x) \bigg[\int u(x,y)\,d\nu(y)\bigg]\,d\nu(x) -
\int f(y) \bigg[\int u(x,y)\,d\nu(x)\bigg]\,d\nu(y) \\
 & = &
\int\!\int (f(x) - f(y))\,u(x,y)\,d\nu(x)\,d\nu(y),
\ene
where in the second last step we changed notations by replacing $x$ with $y$
and $y$ with $x$, using also the property that the product measure
$\nu \otimes \nu$ on $\R^n \times \R^n$ is invariant under the mapping 
$(x,y) \rightarrow (y,x)$. From this we obtain that
$$
\sup_{g \in G} \, \int fg\,d\nu = h \int\!\!\!\int |f(x) - f(y)|\,d\nu(x)\,d\nu(y).
$$

By the Rothaus theorem (cf. \cite{R}, \cite{B-C-G}), the relation \eqref{supfunctional} holds 
true in the class of all smooth bounded $f$ on $\R^n$ if and only if the isoperimetric inequality
$$
\max\{L(1_A),L(-1_A)\} \leq \nu^+(A)
$$
holds true in the class of all Borel (or equivalently, closed) sets $A$ in $\R^n$. Hence, 
the indicator functions play an extremal role in \eqref{supfunctional}. This characterization shows that
the functional forms \eqref{Poincaremunuproductmeas} and \eqref{Poincineqnumean} are equivalent to the same isoperimetric inequality \eqref{isoperimetricmunuproductmeas}.

As a closely related, let us also mention an isoperimetric inequality of Cheeger-type
\begin{equation}
\nu^+(A) \geq h' \min\{\nu(A),1-\nu(A)\},\label{Cheeger}
\end{equation}
in which the optimal value $h'$ is called Cheeger's isoperimetric constant 
associated to the measure $\nu$. A functional form of \eqref{Cheeger} is a Sobolev-type inequality
\begin{equation}
h' \int_{ \R^n} |f - m|\,d\nu \leq \int_{ \R^n} |\nabla f|\,d\nu,\label{Poincineqnumidian}
\end{equation}
where now $m = m(f)$ 
is a median of $f$ under $\nu$, that is, a real number such that
$$
\nu\{f \leq m\} \geq \frac{1}{2}, \quad \nu\{f \geq m\} \geq \frac{1}{2}
$$
(in general the median is not unique, but the left integral in \eqref{Poincineqnumidian} does not depend
on the choice of $m$). The equivalence of \eqref{Cheeger} and \eqref{Poincineqnumidian} is a standard fact.

Comparing \eqref{Cheeger} with \eqref{isoperimetricmunuproductmeas}, it is clear that
$$
h \leq h' \leq 2h.
$$
Often, however, $h' = h$, including the Gaussian measure $\nu = \gamma_n$. 
In this case, by the isoperimetric theorem in Gauss space (\cite{Bor2}, \cite{S-T}), 
the perimeter $\gamma_n^+(A)$ subject to $\gamma_n(A) = \frac{1}{2}$ is 
minimized for any half-space $A_\theta = \{x \in \R^n: \left<x,\theta\right> \leq 0\}$.
But for such sets $\gamma_n^+(A_\theta) = \frac{1}{\sqrt{2\pi}}$, hence 
$h = h' = \sqrt{\frac{2}{\pi}}$. Thus, the constant $\sqrt{\frac{\pi}{2}}$ 
is optimal in the Sobolev-type inequality \eqref{GausPoincp=1}, which may also be stated as the
isoperimetric inequality
$$
\gamma_n^+(A) \geq 2\sqrt{\frac{2}{\pi}}\,\gamma_n(A) (1 - \gamma_n(A))
$$
with an equality for all $A = A_\theta$.

\vskip5mm
\section{{\bf Isoperimetric inequalities for Cauchy measures}}\label{ISOPERIMETRIC INEQUALITIES CAUCHY}
\setcounter{equation}{0}

\vskip2mm
\noindent
Using Lemma \ref{geometricmeaningSobineq}, the Sobolev-type inequality \eqref{PoincareCauchyineqp=1partform} in Corollary \ref{PoincareCauchyp=1} may be equivalently 
stated as an isoperimetric inequality.

\vskip5mm

\begin{corollary}\label{coroisopineqcauc}
Let $\alpha \geq n+1$ and $\beta = \alpha - \frac{n+1}{2}$. 
For any closed set $A$ in $\R^n$,
\begin{equation}
\m_{n,\beta}^+(A) \geq \frac{2\sqrt{\alpha - n}}{\sqrt{\pi}}\,\m_{2n,\alpha}(A \times A^{c}).\label{isopineqcauchy}
\end{equation}
\end{corollary}

\vskip2mm
In order to bring this relation to the form \eqref{isoperimetricmunuproductmeas} as in Lemma \ref{geometricmeaningSobineqproductmeas} and thus prove
Corollary \ref{corisopercauchymeasure}, we need to
bound the Cauchy measure $\m_{2n,\alpha}$ from below in terms of the
product of two Cauchy measures on $\R^n$. To this aim, let us derive
inequality \eqref{boundofcauchy}.

\vskip5mm
\begin{lemma}\label{lemmaboundcauchy}
 For any $\alpha \geq n+1$,
\begin{equation}\label{boundcauchymeas}
\m_{2n,\alpha} \geq d\,\m_{n,\alpha} \otimes \m_{n,\alpha}, \quad
d = d_{n,\alpha} = 
\frac{\Gamma(\alpha - \frac{n}{2})^2}{\Gamma(\alpha - n)\, \Gamma(\alpha)}.
\end{equation}
In particular, $d \geq \frac{1}{2}$ for $\alpha \geq n^2$.
\end{lemma}

\vskip5mm
Thus, if $\alpha$ is sufficiently large, the constant $d$ may be chosen to be
universal.

\vskip2mm
{\bf Proof.} Inequality \eqref{boundcauchymeas} can be stated as a pointwise comparison relation 
for densities of the involved measures. According to the definitions \eqref{cauchyweight} and \eqref{densityCauchyRnRn}, 
the densities of $\m_{2n,\alpha}$ and $\m_{n,\alpha}$ satisfy, for all $x,y \in \R^n$,
\bee
w_{2n,\alpha}(x,y) 
 & = &
\frac{1}{c_{2n,\alpha}}\, \frac{1}{(1 + |x|^2+ |y|^2)^\alpha} \\
 & \geq &
\frac{1}{c_{2n,\alpha}}\, \frac{1}{((1 + |x|^2)(1 + |y|^2))^\alpha} \, = \,
\frac{c_{n,\alpha}^2}{c_{2n,\alpha}}\, w_{n,\alpha}(x) w_{n,\alpha}(y).
\ene
The last product of $w$-functions represents the density
of the product measure $\m_{n,\alpha} \otimes \m_{n,\alpha}$.
In addition, according to \eqref{normalizconst} and \eqref{normalizconstRnRn},
$$
\frac{c_{n,\alpha}^2}{c_{2n,\alpha}} =
\frac{\Gamma(\alpha - \frac{n}{2})^2}{\Gamma(\alpha - n)\, \Gamma(\alpha)}.
$$
This proves \eqref{boundcauchymeas}.

Our next task is to bound the last fraction. To this end,
one may use Stirling's formula. Alternatively, one may appeal to
the following statement proved in \cite{B2}: If
a random variable $\eta>0$ has a log-concave density $q(t)$ on
$(0,\infty)$, then the function $\log \E\,\eta^x - x\log x$
is concave in $x \geq 0$. For example, if $\eta$ has 
a standard exponential distribution with density $q(t) = e^{-t}$, $t>0$, 
we have $\E\,\eta^x = \Gamma(x+1)$, so that the function
$\log \Gamma(x+1) -x\log x$ is concave. Therefore,
\begin{equation}\label{loggamma}
\log \Gamma(x) = \psi(x-1) + u(x), \quad \psi(x) = x\log x, \ \ 
x \geq 1,
\end{equation}
for some concave function $u(x)$. Introduce the operators
$$
\Delta_h U(x) = U(x) - \frac{1}{2}\,U(x-h) - \frac{1}{2}\,U(x+h), \quad h \geq 0.
$$
By Jensen's inequality, if $x-h \geq 1$, then, by \eqref{loggamma},
$$
\Delta_h \log \Gamma(x) \geq \Delta_h \psi(x-1).
$$
In the interval $0 \leq h < x$, the function $\varphi_x(h) = - \Delta_h \psi(x)$ 
has the first two derivatives
\bee
\varphi_x'(h)
 & = &
\frac{1}{2}\,\log(x+h) - \frac{1}{2}\,\log(x-h), \\
\varphi_x''(h) 
 & = &
\frac{1}{2\,(x+h)} + \frac{1}{2\,(x-h)} \, \leq \, \frac{1}{x-h}.
\ene
Since $\varphi_x(0) = \varphi_x'(0) = 0$, Taylor's formula implies
$\varphi_x(h) \leq \frac{h^2}{2(x-h)}$. Hence
$$
\Delta_h \log \Gamma(x) \geq -\varphi_{x-1}(h) \geq - \frac{h^2}{2(x-h-1)}.
$$
Applying this with $x = \alpha - \frac{n}{2}$ and $h = \frac{n}{2}$,
we get
\bee
\log \Gamma\Big(\alpha - \frac{n}{2}\Big) - 
\frac{1}{2}\,\log \Gamma(\alpha - n)  - \frac{1}{2}\,\log \Gamma(\alpha) 
 & \geq &
- \frac{n^2}{8\,(\alpha - n -1)} \\
 & \geq &
- \frac{n^2}{8\,(n^2 - n -1)} \, \geq \, -\frac{1}{2},
\ene
since $\frac{n^2}{n^2 - n - 1} \leq 4$ for $n \geq 2$
(and where we used the assumption $\alpha \geq n^2$). It follows that
\begin{equation}\label{boundratiogamma}
\frac{\Gamma(\alpha - \frac{n}{2})^2}{\Gamma(\alpha - n)\, \Gamma(\alpha)} \geq 
e^{-1/2} > \frac{1}{2}.
\end{equation}

In the case $n=1$, we use the original assumption $\alpha \geq n+1$, implying 
$\alpha \geq 2$. We employ Gautschi's inequality which asserts in particular that
$$
\sqrt{x} < \frac{\Gamma(x + 1)}{\Gamma(x + \frac{1}{2})} < \sqrt{x+1}, \quad
x>0.
$$
Here the left inequality was already used (as Wendel's inequality). Applying
the right inequality with $x = \alpha - 1$, we have
$$
\Gamma\Big(\alpha - \frac{1}{2}\Big) > \frac{1}{\sqrt{\alpha}}\,\Gamma(\alpha).
$$
Hence for the ratio on the left-hand side in \eqref{boundratiogamma} we get
$$
\frac{\Gamma(\alpha - \frac{1}{2})^2}{\Gamma(\alpha - 1)\, \Gamma(\alpha)} >
\frac{\Gamma(\alpha)}{\alpha \Gamma(\alpha - 1)} = \frac{\alpha-1}{\alpha} \geq
\frac{1}{2}.
$$
\qed

\vskip2mm
{\bf Proof of Corollary \ref{corisopercauchymeasure}.}  Note that $\alpha = \beta + \frac{n+1}{2} = \beta^*$, 
so that $\alpha \geq n+1 \Longleftrightarrow \beta \geq \frac{n+1}{2}$, which is
assumed to hold true.

Applying \eqref{boundcauchymeas}, the right-hand side of \eqref{isopineqcauchy} is bounded from below by
$$
\frac{2d}{\sqrt{\pi}}\sqrt{\alpha - n}\ \m_{n,\alpha}(A) \, \m_{n,\alpha}(\bar A).
$$
Since $\alpha - n = \beta - \frac{n-1}{2} \geq \frac{1}{2}\,\beta$,
the desired inequality \eqref{isopineqcauc} follows.

Also, the condition $\alpha \geq n^2$ as in Lemma \ref{lemmaboundcauchy} is fulfilled
as long as $\beta \geq n^2$. In this case,
$$
\frac{2d}{\sqrt{\pi}}\,\sqrt{\alpha - n} \geq 
\frac{1}{\sqrt{\pi}}\sqrt{\beta/2} > \frac{1}{2}\sqrt{\beta}.
$$
\qed

\vskip5mm
\section{\bf Large deviations. Proof of Corollary \ref{corlargedev}}\label{LARGE DEVIATIONS}
\setcounter{equation}{0}

\vskip2mm
\noindent
Let us return to Theorem \ref{Poinccauchytheo} for the general parameter $p$, where it was
assumed that $\alpha > n$ and $1 \leq p < 2\,(\alpha - n)$. It follows from \eqref{Poinccauchy}
that, for any function $f$ on $\R^n$ with Lipschitz semi-norm $\|f\|_{{\rm Lip}} \leq 1$,
\begin{equation}
\int_{\R^n} \int_{\R^n} |f(x) - f(y)|^p\,d\m_{2n,\alpha}(x,y) \, \leq \, 
\frac{1}{\sqrt{\pi}}\,\Big(\frac{\pi}{2}\Big)^p\,
\frac{\Gamma(\frac{p+1}{2})\Gamma(\alpha - n - \frac{p}{2})}{\Gamma(\alpha - n)}.\label{boundbysubunitLipnorm}
\end{equation}
This bound can be used to explore probabilities of large deviations of $f(x) - f(y)$
under the Cauchy measure $\m_{n,\alpha}$. To this aim, it is worthwile to realize how  
the expression on the right-hand side of \eqref{boundbysubunitLipnorm} grows
with respect to the growing parameter $p$. We employ the following
two-sided bound for the Gamma function proposed by Batir \cite{Bat}:
\begin{equation}
\sqrt{2e}\,\Big(\frac{x}{e}\Big)^x \leq \Gamma\Big(x+\frac{1}{2}\Big) \leq
\sqrt{2\pi}\,\Big(\frac{x}{e}\Big)^x, \quad x \geq \frac{1}{2}.\label{Batir}
\end{equation}
Here the constants $\sqrt{2e}$ and $\sqrt{2\pi}$ are optimal for the indicated
$x$-range. As a consequence, whenever $h \geq 0$, $x \geq h + \frac{1}{2}$,
\bee
\frac{\Gamma(x - h + \frac{1}{2})}{\Gamma(x + \frac{1}{2})} 
 & \leq &
\sqrt{\frac{\pi}{e}}\,e^h\,\frac{(x-h)^{x-h}}{x^x} \\
 & = & 
\sqrt{\frac{\pi}{e}}\,e^{h}\,\frac{(1 - \frac{h}{x})^{x-h}}{x^h} \, \leq \,
\sqrt{\frac{\pi}{e}}\,\Big(\frac{e}{x}\Big)^h.
\ene
Applying this with $x = \alpha - n - \frac{1}{2}$ and $h = \frac{p}{2}$
such that $p\leq 2(\alpha-n-1)$, we get
$$
\frac{\Gamma(\alpha - n - \frac{p}{2})}{\Gamma(\alpha - n)} \, \leq \,
\sqrt{\frac{\pi}{e}}\,\Big(\frac{e}{\alpha - n - \frac{1}{2}}\Big)^{p/2}.
$$
Note that since $p \geq 1$, we should strengthen the condition on the range
of $\alpha$ to $\alpha \geq n + \frac{3}{2}$.
In addition, by \eqref{Batir} once more, 
$$
\Gamma\Big(\frac{p+1}{2}\Big) \leq \sqrt{2\pi}\,\Big(\frac{p}{2e}\Big)^{p/2}.
$$
Hence, the right-hand side of \eqref{boundbysubunitLipnorm} can be bounded from above by
$$
\frac{1}{\sqrt{\pi}}\,\Big(\frac{\pi}{2}\Big)^p\,\sqrt{2\pi}\,\Big(\frac{p}{2e}\Big)^{p/2}\,
\sqrt{\frac{\pi}{e}}\,\Big(\frac{e}{\alpha - n - \frac{1}{2}}\Big)^{p/2} <
2\,\Big(\frac{\frac{\pi^2}{8}\,p}{\alpha - n - \frac{1}{2}}\Big)^{p/2} \leq
2\,\Big(\frac{2 p}{\alpha - n}\Big)^{p/2},
$$
where in the last step we used that $\alpha - n - \frac{1}{2} \geq \frac{2}{3}\,(\alpha - n)$ and
$\frac{3\pi^2}{16} < 2$. 
One can summarize.

\vskip5mm
\begin{corollary}\label{corboundbysubunitLipnorm}
 Let $\alpha\geq n+\frac{3}{2}$ and $1\leq p\leq 2(\alpha-n-1)$. 
For any function $f$ on $\R^n$ with $\|f\|_{\rm Lip} \leq 1$,
\begin{equation} \label{boundbysubunitLipnorm2}
\int_{\R^n} \int_{\R^n} |f(x) - f(y)|^p\,d\m_{2n,\alpha}(x,y) \, \leq \, 
2\,\Big(\frac{2p}{\alpha - n}\Big)^{p/2}.
\end{equation}
\end{corollary}

\vskip2mm
{\bf Proof of Corollary \ref{corlargedev}.} By Markov's inequality, for any $t>0$ and the range
of $(n,p)$ as in Corollary \ref{corboundbysubunitLipnorm}, from \eqref{boundbysubunitLipnorm2} we get
$$
\m_{2n,\alpha}\big\{\sqrt{\alpha - n}\, |f(x) - f(y)| \geq t\big\} \leq
2\,\Big(\frac{2p}{t^2}\Big)^{p/2}.
$$
Choosing here $p = t^2/(2e)$, the right-hand side becomes
$
2\,e^{-p/2} = 2\,e^{-t^2/(4e)} \leq 2\,e^{-t^2/12},
$
which is applicable for $2e \leq t^2 < 4e\,(\alpha - n - 1)$.
Thus, for this interval
\begin{equation}\label{subgaussestimate}
\m_{2n,\alpha}\Big\{\sqrt{\alpha - n}\, |f(x) - f(y)| \geq t\Big\} \leq
2\,e^{-t^2/12}.
\end{equation}
This inequality continues to hold for $t < t_0 = \sqrt{2e}$, since
$2\,e^{-t_0^2/12} > 1$, while the left probability in \eqref{subgaussestimate} does not
exceed 1. lso, the right end of the interval may be simplified using
$\alpha - n - 1 \geq \frac{1}{3}\,(\alpha - n)$ and $\frac{2}{3}\sqrt{e} > 1$. 
As a result, the subgaussian deviation inequality \eqref{subgaussestimate}, that is, \eqref{subgausineq} 
holds true in the interval $0 \leq t \leq \sqrt{\alpha - n}$.

\noindent If additionally $\alpha \geq n^2$, then, by Lemma \ref{lemmaboundcauchy}, 
$\m_{2n,\alpha} \geq \frac{1}{2}\,\m_{n,\alpha} \otimes \m_{n,\alpha}$,
which gives \eqref{subgausineq2}.
\qed

\vskip5mm
\section{{\bf Convexity of Cauchy measures}}\label{CONVEXITY}
\setcounter{equation}{0}

\vskip2mm
\noindent
The class of Cauchy measures have been intensively studied in the literature
in the framework of convex (hyperbolic) measures. Let us recall several
results in this direction.
 
A probability measure $\nu$ on $\R^n$ is called $\kappa$-concave, where 
$-\infty \leq \kappa \leq \infty$, if it satisfies the Brunn-Minkowski-type inequality
\begin{equation}\label{BrunnMink}
\nu(tA + (1 - t)B) \geq \big(t\nu(A)^\kappa + (1 - t) \nu(B)^\kappa\big)^{1/\kappa}
\end{equation}
for all $t \in (0,1)$ and for all Borel measurable sets $A,B \subset \R^n$ with positive measure.
Here $tA + (1 - t)B = \{tx + (1 - t)y : x \in A, y \in B\}$ stands for the Minkowski sum
of the two sets.

Since the right-hand side of \eqref{BrunnMink} is non-decreasing in $\kappa$,
The class of $\kappa$-concave measure is getting smaller for the growing
parameter $\kappa$ and contains only $\delta$-measure in the limit case
$\kappa = \infty$. Otherwise, necessarily $\kappa \leq 1$.
The Lebesgue measure and its restriction to convex bodies in $\R^n$
are $\kappa$-concave with $\kappa = \frac{1}{n}$. This is the content of the
Brunn-Minkowski-Lyusternik theorem. If $\kappa = 0$, we obtain
the important class of log-concave measures, in which case \eqref{BrunnMink} takes the form
$$
\nu(tA + (1 - t)B) \geq \nu(A)^t \nu(B)^{1-t}.
$$
The largest class corresponding to the other limit value $\kappa = -\infty$ is described by
the Brunn-Minkowski-type inequality
$$
\nu(tA + (1 - t)B) \geq \min\{\nu(A),\nu(B)\}.
$$
Such measures are called convex (following Borell) or hyperbolic (according to V. Milman).

A characterization of $\kappa$-concave measures was given by Borell \cite{Bor1}.
For simplicity suppose that  $\nu$ is absolutely continuous with respect to the Lebesgue 
measure on $\R^n$. Then any convex measure should be supported on an open convex set 
$\Omega \subset \R^n$ where it has a continuous density. Moreover,
the log-concavity of $\nu$ is equivalent to the property that its density
is log-concave on $\Omega$. In the case $\kappa<0$, $\nu$ is $\kappa$-concave 
if and only if it has density of the form
$$
w(x) = \frac{d\nu(x)}{dx} = V(x)^{-p}, \quad \kappa = - \frac{1}{p - n} \ (p \geq n),
$$
for some convex function $V$ on $\Omega$.

The $n$-dimensional Cauchy measure $\nu = \m_{n,\alpha}$, $\alpha>\frac{n}{2}$,
corresponds to the convex function $V(x) = \sqrt{1 + |x|^2}$. Hence this measure
is $\kappa$-concave with the optimal parameter 
$$
\kappa = - \frac{1}{2\alpha - n}.
$$ 

It was shown in \cite{B3} that any $\kappa$-concave probability measure $\nu$ on $\R^n$
with $\kappa \leq 1$ satisfies an isoperimetric inequality
\begin{equation}
\nu^+(A) \geq \frac{c(\kappa)}{m}\,\big(\nu(A) (1 - \nu(A))\big)^{1 - \kappa},\label{isopineq2}
\end{equation}
where $m$ is the $\nu$-median of the Euclidean norm $x \rightarrow |x|$, and where 
$c(\kappa)$ is a positive continuous function in the range $(-\infty,1]$.
A closely related weighted Sobolev-type inequality
\begin{equation}\label{weightedSobolev}
\int_{\R^n} \int_{\R^n} |f(x) - f(y)|\,d\nu(x)\,d\nu(y) \, \leq \, C(\kappa)
\int_{\R^n} |\nabla f(x)| \, (m - \kappa x)\, d\nu(x)
\end{equation}
was derived in \cite{B-L2} for $\kappa \leq 0$, where $m = \exp \int \log |x|\,d\nu(x)$
is the geometric mean (or $L^0$-norm) of the Euclidean norm with respect to $\nu$.

If $\nu$ is log-concave, these inequalities are equivalent within universal
factors, and \eqref{weightedSobolev} is reduced to the non-weighted $L^1$-Poincar\'e-type inequality
\begin{equation}
\int_{\R^n} \int_{\R^n} |f(x) - f(y)|\,d\nu(x)\,d\nu(y) \, \leq \, Cm
\int_{\R^n} |\nabla f(x)| \,d\nu(x).\label{L1Poincagain}
\end{equation}
When $\nu$ is a uniform distribution on a convex body, \eqref{L1Poincagain}
corresponds to the result by Kannan, Lov\'asz and Simonovits \cite{K-L-S};
the general log-concave case $\kappa = 0$ was considered in \cite{B2}. 
Here the quantity $m$ is equivalent to the $L^1$-norm
$\int |x|\,d\mu(x)$. Hence, when $\nu$ is isotropic with
$\int |x|^2\,d\textcolor{red}{\nu}(x) = 1$, (11.4) yields 
$$
\int_{\R^n} \int_{\R^n} |f(x) - f(y)|\,d\nu(x)\,d\nu(y) \, \leq \, C\sqrt{n}
\int_{\R^n} |\nabla f(x)| \,d\nu(x).
$$
This is far from beeing optimal, since the factor
$\sqrt{n}$ can be replaced with $\sqrt{\log n}$, due to a recent result of
Klartag \cite{K}.

As for the general case $\kappa < 0$, the inequalities \eqref{isopineq2}-\eqref{weightedSobolev} cannot be used
to recover the concentration of measure phenomenon, such as the one
for the Gaussian measure. In this respect, known results about Cauchy measures
are more accurate. For example. it was shown in \cite{B-L1} that every $\m_{n,\alpha}$, 
$\alpha \geq n$, admits a weighted Poincar\'e-type inequality
\begin{equation}\label{weightedPoincar}
\Var_{\m_{n,\alpha}}(f) \leq \frac{C_\alpha}{2(\alpha - 1)} 
\int |\nabla f(x)|^2\,(1 + |x|^2)\,d\m_{n,\alpha}(x)
\end{equation}
in the class of all smooth functions $f$ on $\R^n$ with constant
$C_\alpha = \big(\sqrt{1 + \frac{2}{\alpha - 1}} + \sqrt{\frac{2}{\alpha - 1}}\,\big)^2$ (see also \cite{BBDGV} \cite[Theorem 2]{BDGV} for the optimal values of $C_{\alpha}$ in the whole range of admissible $\alpha$'s).
Since $C_\alpha \rightarrow 1$ as $\alpha \rightarrow \infty$,
after the linear rescaling of the space variable, \eqref{weightedPoincar} yields in the limit
the Gaussian Poincar\'e-type inequality
$$
\Var_{\gamma_n}(f) \leq \int |\nabla f(x)|^2\,d\gamma_n(x).
$$
See also \cite{BDGV}
As was already mentioned, Theorem \ref{Poinccauchytheo} with $p=1$ also yields
the $L^1$-Poincar\'e-type inequality for the Gaussian measure with a dimension free
optimal constant in the limit as $\alpha \rightarrow \infty$.

\vskip4mm
{\bf Acknowledgements.} This work was started in 2018 when B.V. visited to the
University of Minnesota and was continued in 2022 when S.B. visited to the
Parthenope University of Naples. The authors are grateful for hospitality.
Research of the first author was partially supported by the NSF grant DMS-2154001.
B.V. was partially supported by Gruppo Nazionale per l'Analisi Matematica, la Probabilit\`a e le loro Applicazioni (GNAMPA) of Istituto Nazionale di Alta Matematica (INdAM).
This study was also carried out within the ``Geometric-Analytic Methods for PDEs and Applications (GAMPA)'' projects - funded by the Ministero dell'Universit\`a e della Ricerca - within the PRIN 2022 program (D.D.104 - 02/02/2022). This manuscript reflects only the authors' views and opinions and the Ministry cannot be considered responsible for them.

\end{document}